\documentclass[11pt,a4paper,reqno]{amsart}
\numberwithin{equation}{subsection}
\usepackage[utf8]{inputenc}
\usepackage[english]{babel}
\usepackage{hyperref}

\usepackage{multirow} 
\usepackage{graphicx}
\usepackage{array}
\usepackage{longtable}
\usepackage{tikz-cd}
\DeclareMathOperator{\lcm}{lcm}

\def\ZZ{\mathbb{Z}}

\usepackage[foot]{amsaddr}
\usepackage{url}
\usepackage{hyperref}

\usepackage{amsmath}
\usepackage{amssymb, amscd}
\usepackage[all, cmtip]{xy}
\usepackage{xcolor}
\usepackage{multirow} 
\usepackage{longtable}
\usepackage{array}

\newtheorem{theorem}{Theorem}[section]

\newtheorem{conjecture}[theorem]{Conjecture}

\newtheorem{corollary}[theorem]{Corollary}
\newtheorem{lemma}[theorem]{Lemma}

\theoremstyle{definition}
\newtheorem{definition}[theorem]{Definition}

{\sc}
\newtheorem{example}[theorem]{Example}
{}
{}
{\sc}
\newtheorem{remark}{Remark}[section]

\def\ZZ{\mathbb{Z}}

\newcommand{\thismonth}{\ifcase\month\or
  January\or February\or March\or April\or May\or June\or
  July\or August\or September\or October\or November\or December\fi
  \space\number\year}

\makeatletter
\newcommand{\rssymb}[2]{\newcommand{#1}{{\mathrmsl{#2}}}}
\newcommand{\calsymb}[2]{\newcommand{#1}{{\mathcal{#2}}}}
\newcommand{\bbsymb}[2]{\newcommand{#1}{{\mathbb{#2}}}}
\newcommand{\lieoper}[2]{\newcommand{#1}{\mathop
  {\mathfrak{#2}\null}\nolimits}}
\newcommand{\oper}[3][n]{\newcommand{#2}{\mathop
  {\mathrm{#3}\null}\ifx n#1\nolimits\else\limits\fi}}
\newcommand{\rsoper}[3][n]{\newcommand{#2}{\mathop
  {\mathrmsl{#3}\null}\ifx n#1\nolimits\else\limits\fi}}
\bbsymb\C{C} \bbsymb\F{F} \bbsymb\HQ{H}\bbsymb\N{N} \bbsymb\Q{Q}
\bbsymb\R{R} \bbsymb\U{U} \bbsymb\V{V} \bbsymb\W{W} \bbsymb\Z{Z}
\bbsymb\bbf{F} \bbsymb\bbk{K} \bbsymb\bbi{I} \bbsymb\bbl{L}
\bbsymb\bbo{O} \bbsymb\bbj{J} \bbsymb\bby{Y} \bbsymb\bbp{P}
\bbsymb\bba{A}
\calsymb\cA{A} \calsymb\cB{B} \calsymb\cC{C} 
\calsymb\cM{M} \calsymb\cN{N} \calsymb\cO{O} \calsymb\cP{P}
\calsymb\cU{U} \calsymb\cV{V} \calsymb\cW{W} \calsymb\cX{X}
\calsymb\cY{Y} \calsymb\cZ{Z}
\renewcommand{\geq}{\geqslant} \renewcommand{\leq}{\leqslant}
\oper\End{End} \oper\Hom{Hom}                    
\oper\Sym{Sym} \oper\Skew{Skew}
\oper\Aut{Aut}                                   
\oper\GL{GL} \oper\SL{SL}\oper\Symp{Sp} \oper\CO{CO} \oper\On{O}
\oper\SO{SO} \oper\Pin{Pin} \oper\Spin{Spin} \oper\CU{CU}
\oper\Un{U} \oper\SU{SU} \oper\PSU{PSU} \rsoper\Diff{Diff}
\rsoper\SDiff{SDiff}
\lieoper\der{der}                                
\lieoper\gl{gl} \lieoper\sgl{sl}\lieoper\symp{sp} \lieoper\co{co}
\lieoper\so{so} \lieoper\spin{spin} \lieoper\cu{cu} \lieoper\un{u}
\lieoper\su{su} \rsoper\Vect{Vect} \rsoper\Ham{Ham}
\def\la#1{\hbox to #1pc{\leftarrowfill}}
\def\ra#1{\hbox to #1pc{\rightarrowfill}}

\newcommand{\Norm}[2][]{\bigl|\mkern-3mu\bigr|#2\bigr|\mkern-3mu\bigr|
  _{\lower1pt\hbox{${}_{#1}$}}}

\rsoper\dimn{dim}                           
\rsoper\grad{grad}                          
\rsoper\kernel{ker}\rsoper\image{im}        
\rsoper\alt{alt}   \rsoper\sym{sym}         
\rsoper\Ad{Ad}     \rsoper\ad{ad}           
\rsoper\CoAd{CoAd} \rsoper\coad{coad}       
\rsoper\trace{tr}  \rsoper\trfree{tf}       
\rsoper\detm{det}                           
\rsoper\Vol{Vol}                            
\rsoper\divg{div}                           
\rsoper\sign{sign}                          
\rssymb\iden{id}                            
\rssymb\vol{vol}                            
\oper\Imag{Im}\oper\Real{Re}                
\newcommand{\sd}{{\raise1pt\hbox{$\scriptscriptstyle +$}}}
\newcommand{\asd}{{\raise1pt\hbox{$\scriptscriptstyle -$}}}
\newcommand{\sdasd}{{\raise1pt\hbox{$\scriptscriptstyle\pm$}}}

\newcommand{\asdsd}{{\raise1pt\hbox{$\scriptscriptstyle\mp$}}}

\rsoper\scal{scal}
\def\kahl/{k\"ahler}
\def\Kahl/{K{\"a}hler}

\begin{document}

\title[Berglund-Hübsch transpose and Sasaki-Einstein  Rational homology 7-spheres]
{Berglund-Hübsch transpose and Sasaki-Einstein  Rational homology 7-spheres}

\author[J. Cuadros]{Jaime Cuadros Valle$^1$}
\author[R.Gomez]{Ralph R. Gomez$^2$}
\author[J. Lope]{Joe Lope Vicente$^1$}

\address{$^1$Departamento de Ciencias, Secci\'on Matem\'aticas,
Pontificia Universidad Cat\'olica del Per\'u,
Apartado 1761, Lima 100, Per\'u}
\address{$^2$Department of Mathematics and Statistics, Swarthmore College, 500 College Avenue, Swarthmore, PA 19081, USA}
\email{jcuadros@pucp.edu.pe}
\email{rgomez1@swarthmore.edu}
\email{j.lope@pucp.edu.pe}

\date{\thismonth}

\begin{abstract} 
We show that links of invertible polynomials coming from the Johnson and  Koll\'ar list of K\"ahler-Einstein 3-folds that are  rational homology 7-spheres remain rational homology 7-spheres under the so-called Berglund-Hübsch transpose rule coming from classical mirror symmetry constructions. Actually, this rule produces twins, that is, links with same degree,  Milnor number and homology $H_3$, with the exception of iterated Thom-Sebastiani sums of singularities of chain and cycle type, where the torsion and the Milnor number vary. The Berglund-Hübsch  transpose rule not only gives a framework to better understand the existence of Sasaki-Einstein twins but also gives a mechanism for producing new examples of Sasaki-Einstein twins in the rational homology 7-sphere setting. We  also give reasonable conditions for a Sasaki-Einstein rational homology 7-sphere  to remain Sasaki-Einstein under the Berglund-Hübsch transpose rule. In particular, we found 75 new examples of Sasaki-Einstein rational homology 7-spheres arising as links of not well-formed hypersurface singularities.

\end{abstract}

\maketitle

\medskip

\medskip

\maketitle
\vspace{-2mm}


\section{introduction}

Recall that an odd dimensional compact Riemannian manifold is Sasaki if we have the quadruple of tensors $(\phi, \xi, \eta, g)$ where $\eta$ is a contact 1-form, $\xi$ is the Reeb vector field, $\phi$ is an endomorphism of the tangent bundle satisfying $\phi^{2}=-I+\xi\otimes \eta $ and $g$ is a Riemannian metric compatible with the structure tensors. Moreover, $\xi$ is a Killing vector field and when $\phi$ is restricted to the contact sub-bundle $D=\ker \eta$, it is integrable,  that is, the underlying CR-structure is integrable. We say $M$ is Sasaki-Einstein if in addition, the Riemannian metric $g$ is Einstein, that is, $\mathrm{Ric}_g=\lambda g.$ Sasaki geometry has been extremely useful in constructing a multitude of odd dimensional Riemannian manifolds which are Einstein of positive scalar curvature. 
Besides its inherent importance in differential geometry, Sasaki-Einstein metrics are of great interest in mathematical physics, since they admit real Killing spinors which provide many applications in theoretical physics. 
For instance,  they  play an important role in the context of superstring theory and  $M$-theory \cite{FK}.  Also a string theory conjecture known as the AdS/CFT correspondence relates superconformal field theories and Sasaki-Einstein metrics in dimension five and seven \cite{GMS}, \cite{XY}.  These metrics are abundant (see  \cite{BBG} and references therein): in  dimension 5, it is known that there exist Sasaki-Einstein structures on $\# k\left(S^2 \times S^3\right)$ for all values of $k$.  Moreover, $S^5$  admits many inequivalent Sasaki-Einstein structures.  There are also quasi-regular Sasaki-Einstein metrics on 5-manifolds which are not connected sums of $S^2 \times S^3$, which include infinitely many rational homology 5-spheres. In dimension 7, all 28 oriented diffeomorphism classes on $S^7$ admit Sasaki-Einstein metrics. In \cite{BGN2, CL} the existence of  Sasaki-Einstein metrics in rational homology 7-spheres were established. Sasaki-Einstein metrics on connected sums of $\#k(S^3\times S^4)$ were found for 22 different values of $k$ in \cite{CL} which, in particular,  provide examples of linearly unstable Sasaki-Einstein 7-manifolds since they admit non-trivial harmonic 3-forms, see \cite{SWW} for details.\\ 
\indent One reason for the plethora of examples of Sasaki-Einstein manifolds is due to the construction developed in \cite{BG1} involving links of an isolated hypersurface singularity at the origin. In particular,  hypersurfaces defined by invertible polynomials have proven to be quite fruitful when used with this technique. These are weighted homogeneous polynomials $f \in \mathbb{C}\left[x_0, \ldots, x_n\right]$, which are a sum of exactly $n$ monomials, such that the weights $w_1, \ldots, w_n$ of the variables $x_1, \ldots, x_n$ are unique  and that  its affine cone is smooth outside  $(0, \ldots, 0),$ {\it i.e.,} $f$ is quasismooth.  The class of invertible polynomials includes all polynomials of Brieskorn-Pham type, but  do not refer only to them. These polynomials were  studied by Berglund and Hübsch \cite{BH}, who showed that a mirror $n$-fold is related to a dual polynomial $f^T$. For an invertible polynomial $f=\sum_{j=1}^n \prod_{i=1}^n x_i^{a_{i j}}$ the transpose or dual  polynomial $f^T$ is defined by transposing the exponential  matrix $A=\left(a_{i j}\right)_{i, j}$ of the original polynomial, thus $f^T=\sum_{j=1}^n \prod_{i=1}^n x_i^{a_{j i}}$. Notice that if the polynomial is of Brieskorn-Pham type then the polynomial is self-dual. The idea of how to construct these mirror pairs was made precise by Krawitz et al., the interesting reader is referred to \cite{Kr, KPABR}.

In this paper we continue the study of  Sasaki-Einstein structures on 7-manifolds through the 
study of the $S^1$ orbibundles over the 1936 Fano 3-folds found by Johnson and Kollár \cite{JK} that admit Kähler-Einstein metrics, as done  initially  in  \cite{BGN2}  and continued in \cite{CL}. 
For most of the elements given in list in \cite{JK},  for the given data $(d, \mathbf{w}=\left(w_{0}, \ldots, w_{4}\right))$, one can find different types of invertible polynomials  representing them. 
It is natural to apply the Berglund-Hübsch transpose rule to each possible polynomial representations of the given data as done previously in \cite{Go}. We focus on data $(d, \mathbf{w})$  that yield to  links that are  Sasaki-Einstein rational homology 7-spheres. Under these circumstances, we find that   Berglund-Hübsch  transpose method gives rise to rational homology 7-spheres. Furthermore,  we prove  that  via this method, new rational homology 7-spheres outside the list given in \cite{BGN2} and \cite{CL}  can only be produced in case the polynomial representation is of chain-cycle type. Finally, we give weak conditions for this new links to admit Sasaki-Einstein metrics. Actually, thanks to  Theorem 4.2, we established the existence of  75 new rational homology 7-spheres admitting Sasaki-Einstein metrics that come from not well-formed orbifolds. We also notice that if one begins with a Fano orbifold, for whom existence of K\"ahler-Einstein metric is uncertain, the output via the Berglund-Hübsch transpose rule may produce Sasaki-Einstein rational homology spheres, see Example 4.5. In addition, we find very mild conditions, such that the Berglund-H\"ubsh transpose produces {\it twins}, that is, $2$-connected rational homology 7-spheres  with identical homologies, degrees and  Milnor numbers admitting Sasaki-Einstein metrics. 

The paper is organized as follows: in Section 2 we present background material which includes some basics on Sasaki-Einstein geometry and its natural presence on certain links of hypersurface singularities and discuss how to calculate the homology of the links; we also briefly describe  the Berglund-Hubsch transpose rule from mirror symmetry and its application to Sasakian geometry. In Section 3 and Section 4 we state and prove our main results. In section 5  we include a table listing  new examples of Sasaki-Einstein rational homology 7-spheres provided by Theorem 4.2,   and links to codes written in Matlab to ease  computations.  

\section{preliminaries}

\subsection{Links and Sasakian geometry} Let us consider the weighted $\mathbb{C}^{*}$-action on $\mathbb{C}^{n+1}$ given by 
$$
\left(z_{0}, \ldots, z_{n}\right) \longmapsto\left(\lambda^{w_{0}} z_{0}, \ldots, \lambda^{w_{n}} z_{n}\right)
$$
where $w_{i}$ are the weights which are positive integers and $\lambda \in \mathbb{C}^{*}$. Let  $\mathbf{w}=\left(w_{0}, \ldots, w_{n}\right)$ be the weight vector such that 
$
\operatorname{gcd}\left(w_{0}, \ldots, w_{n}\right)=1.
$
Recall that a weighted homogeneous polynomial  $f \in \mathbb{C}\left[z_{0}, \ldots, z_{n}\right]$  of degree $d$ and weight $\mathbf{w}=\left(w_{0}, \ldots, w_{n}\right)$ satisfies that for any $\lambda \in \mathbb{C}^{*}=\mathbb{C} \backslash\{0\}$,
$$
f\left(\lambda^{w_{0}} z_{0}, \ldots, \lambda^{w_{n}} z_{n}\right)=\lambda^{d} f\left(z_{0}, \ldots, z_{n}\right) .
$$
We will assume that  $f$ is quasi-smooth, that is, its zero locus in $\mathbb{C}^{n+1}$ has only an isolated singularity at the origin. 
The link $L_{f}(\mathbf{w}, d),$ or  $L_f$ for short,  is defined by  $$L_{f}(\mathbf{w}, d)=f^{-1}(0) \cap S^{2 n+1},$$ where $S^{2 n+1}$ is the $(2 n+1)$-sphere in $\mathbb{C}^{n+1}$. It follows from the Milnor fibration theorem \cite{Mi}  that $L_{f}(\mathbf{w}, d)$ is a closed $(n-2)$-connected $(2n-1)$-manifold that bounds a parallelizable manifold with the homotopy type of a bouquet of $n$-spheres. It is well-known \cite{BBG} that  $L_{f}(\mathbf{w}, d)$ admits a quasi-regular Sasaki structure  $\mathcal{S}_{\mathbf{w}}=\left(\xi_{\mathbf{w}}, \eta_{\mathbf{w}}, \Phi_{\mathbf{w}}\right)$, which is the restriction of the weighted Sasakian on the sphere $S^{2 n+1}$ with Reeb vector field $\xi_{\mathrm{w}}=\sum_k w_k H_k$ where $H_k=-i\left(z_k \partial_{z_k}-\right.$ $\left.\bar{z}_k \partial_{\bar{z}_k}\right)$. Moreover, if one considers the locally free $S^1$-action induced by the weighted $\mathbb{C}^{*}$ action on $f^{-1}(0)$ the quotient space of the link  $L_{f}(\mathbf{w}, d)$ by this action is the weighted hypersurface $X_{f},$ a K\"ahler orbifold: we have the following commutative diagram \cite{BBG}

\begin{equation*}
\begin{CD} 
 L_{f}(\mathbf{w}, d) @> {\qquad\qquad}>>  S^{2n+1}_{\bf w}\\
@VV{\pi}V  @VVV\\
X_{f}  @> {\qquad\qquad}>>  {\mathbb P}({\bf w}),
\end{CD}
\end{equation*}
where $S_{\mathbf{w}}^{2 n+1}$ denotes the unit sphere with a weighted Sasakian structure, $\mathbb{P}(\mathbf{w})$ is a weighted projective space coming from the quotient of $S_{\mathbf{w}}^{2 n+1}$ by a weighted circle action generated from the weighted Sasakian structure. The top horizontal arrow is a Sasakian embedding and the bottom arrow is a Kählerian embedding and  the vertical arrows are orbifold Riemannian submersions. 

It follows from the orbifold adjunction formula that the quotient orbifold $X_{f}$ by the natural  $S^1$-action is Fano if  
\begin{equation*}
I=|{\bf w}|-d_f >0.  
\end{equation*}
Here $|\mathbf{w}|=\sum_{i=0}^m w_i$ denotes the norm of the weight vector $\mathbf{w}$ and $d_f$ is the degree of the polynomial $f$, one refers to its difference $I$ as the index.  

In \cite{Ko}, Kobayashi showed that  the link of a cone over a smooth projective variety $X \subset \mathbb{P}^N$ carries a natural Einstein metric if and only if $X$ is Fano and $X$ carries a Kähler-Einstein metric. In \cite{BG1}, the authors generalized this result to weighted cones and furthermore gave an algorithm, the Kobayashi-Boyer-Galicki method, to obtain $(n-1)$-connected Sasaki-Einstein $(2n+1)$-manifolds from the existence of  orbifold Fano Kähler-Einstein hypersurfaces $X_{f}$ in weighted projective $2n$-space $\mathbb{P}(\mathbf{w}).$  To be more precise, we  have the following result.

\begin{theorem}[\cite{Bo}]Let $L_f$ be a link of an isolated hypersurface singularity defined by a quasi-homogeneous polynomial $f$ of degree $d$ with a Sasakian structure so that we have the $S^1$ orbibundle $S^1 \rightarrow L_f \rightarrow X_f$. If $X_f$ is a Fano hypersurface of degree $d$  in weighted projective space $\mathbb P(\mathbf{w})$ with $I=|\mathbf{w}|-d>0$, then $L_f$ admits a Sasakian structure such that the Sasakian metric has positive Ricci curvature. Moreover, if $X_f$ is a Fano orbifold hypersurface in $\mathbb{P}(\mathbf{w})$ satisfying the Sasaki-Einstein inequality
$$
I d<\frac{n}{n-1} \min _{i, j}\left(w_i w_j\right),
$$
then $L_f$ admits a Sasaki-Einstein metric.
\end{theorem}
\medskip

Thus, links of isolated hypersurface singularities provide a method for constructing Sasaki-Einstein metrics on certain odd dimensional manifolds via the existence of special hypersurfaces in weighted projective space. 
In \cite{JK}, Johnson and Kollár give a list of sporadic 4442 quasi-smooth Fano 3 -folds $\mathcal{Z}$ anticanonically embedded in weighted projective 4-spaces $\mathbb{P}(\mathbf{w})$. Moreover, they show that 1936 of these $3$-folds admit Kähler-Einstein metrics. 
In \cite{BGN2} and  \cite{CL}  links associated to this K\"ahler-Einstein $3$-folds were studied. Applying Orlik's formula, the topology of these links were computed. From that lot, there are 236  Sasaki-Einstein rational homology 7-spheres and all but one of these quasismooth hypersurfaces  of degree $d$ embedded in $\mathbb{P}(w_0, w_1, w_2, w_3,w_4)$  with $w_0\leq w_1 \leq w_2 \leq w_3 \leq w_4,$  satisfy the inequality
\begin{equation}
d \leq w_0 w_1.
\end{equation}
The divergent element is the 3-fold defined by the weight vector (13,143,775,620,465) and degree $d=2015,$ (the first element in Table 1 in \cite{CL}).
It is known (see Proposition 3.3 in \cite{JK}) that if the inequality (2.1.1) is satisfied, then the quasismooth hypersurface 
$X\subset\mathbb{P}(w_0, w_1, w_2, w_3,w_4)$ does not admit a {\it tiger}. (Recall, a tiger is an effective $\mathbb{Q}$-divisor $D$ such that $D$ is numerically equivalent to the anticanonical divisor $-K_X$ and $(X, D)$ is not Kawamata log-terminal.) 

\subsection{Topology of the links}  Recall that 
the Alexander  polynomial $\Delta_f(t)$ in \cite{Mi} associated to a link $L_f$ of dimension $2 n-1$ is the characteristic polynomial of the monodromy map $$h_*: H_{n}(F, \mathbb{Z}) \rightarrow H_{n}(F, \mathbb{Z})$$  induced by the circle action on the Milnor fibre $F$. Then, 
$\Delta_f(t)=\operatorname{det}\left(t {\mathbb I}-h_*\right)$. 
Now both $F$ and its closure $\bar{F}$ are homotopy equivalent to a bouquet of $n$-spheres $S^n \vee \cdots \vee S^n,$ and the boundary of $\bar{F}$ is the link $L_f$, which is $(n-2)$-connected. The Betti numbers $b_{n-1}\left(L_f\right)=b_{n}\left(L_f\right)$ equal the number of factors of $(t-1)$ in $\Delta_f(t)$.  From the  Wang sequence of the Milnor fibration (see \cite{Or}) $$\quad 0 \longrightarrow H_n\left(L_f, \mathbb{Z}\right) \longrightarrow H_n(F, \mathbb{Z}) \stackrel{\mathbb{I}-h_*}{\longrightarrow} H_n(F, \mathbb{Z}) \longrightarrow H_{n-1}\left(L_f, \mathbb{Z}\right) \longrightarrow 0$$ one can derive the following facts:
\begin{enumerate}
\item $L_f$ is a rational homology sphere if and only if $\Delta_f(1) \neq 0.$
\item $L_f$  homotopy sphere if and only if $\Delta_f(1)= \pm 1.$ 
\item If  $L_f$ is a rational homology sphere, then the order of $H_{n-1}\left(L_f, \mathbb{Z}\right)$ equals $|\Delta_f(1)|.$ 
\end{enumerate}


In the case that $f$ is a weighted homogeneous polynomial  there is an algorithm due to Milnor and Orlik \cite{MO} to calculate  the free part of $H_{n-1}\left(L_f, \mathbb{Z}\right).$ The authors associate to any monic polynomial $f$ with roots $\alpha_1, \ldots, \alpha_k \in \mathbb{C}^*$ its divisor
$$
\operatorname{div} f=\left\langle\alpha_1\right\rangle+\cdots+\left\langle\alpha_k\right\rangle
$$
as an element of the integral ring $\mathbb{Z}\left[\mathbb{C}^*\right].$ Let  $\Lambda_n=\operatorname{div}\left(t^n-1\right)$. 
Then the divisor of $\Delta_f(t)$ is given by
\begin{equation}
\operatorname{div} \Delta_f=\prod_{i=0}^n\left(\frac{\Lambda_{u_i}}{v_i}-\Lambda_1\right),
\end{equation}
where the $u_i's$  and $v_i's$ are given  terms of the degree $d$ of $f$ and the weight vector ${\bf w}=(w_0,\ldots w_n)$ by the equations 
\begin{equation}
u_i=\frac{d}{\operatorname{gcd}\left(d, w_i\right)}, \quad v_i=\frac{w_i}{\operatorname{gcd}\left(d, w_i\right)}.
\end{equation}

Using the relations $\Lambda_a \Lambda_b=\operatorname{gcd}(a, b) \Lambda_{\operatorname{lcm}(a, b)}$, Equation (2) takes the form
\begin{equation}
\operatorname{div} \Delta_f=(-1)^{n+1}\Lambda_1+\sum a_j \Lambda_j,
\end{equation}

where $a_j \in \mathbb{Z}$ 
 and the sum is taken over the set of all least common multiples of all combinations of the $u_0, \ldots, u_n$. Then the Alexander polynomial has an alternative expression given by
$$
\Delta_f(t)=(t-1)^{(-1)^n} \prod_j\left(t^j-1\right)^{a_j},
$$
and

\begin{equation}
b_{n-1}\left(L_f\right)=(-1)^{n+1}+\sum_j a_j.
\end{equation}

Also is well-known (see \cite{MO}) that the Milnor number is given by $$\left(\frac{d-w_0}{w_0}\right)\left(\frac{d-w_1}{w_1}\right)\left(\frac{d-w_2}{w_2}\right)\left(\frac{d-w_3}{w_3}\right)\left(\frac{d-w_4}{w_4}\right)=\mu.$$

Moreover, Milnor and Orlik gave and explicit formula to calculate the free part of $H_{n-1}\left(L_f, \mathbb{Z}\right):$ 
\begin{equation}
b_{n-1}\left(L_{f}\right)=\sum(-1)^{n+1-s} \frac{u_{i_{1}} \cdots u_{i_{s}}}{v_{i_{1}} \cdots v_{i_{s}} \operatorname{lcm}\left(u_{i_{1}}, \ldots, u_{i_{s}}\right)},
\end{equation}
where the sum is taken over all the $2^{n+1}$ subsets $\left\{i_{1}, \ldots, i_{s}\right\}$ of $\{0, \ldots, n\}$. In \cite{Or}, Orlik gave a conjecture (\cite{Or}) which allows to determine the torsion of the homology groups  of the link  in terms of the weight of $f.$

\begin{conjecture}[Orlik] 
Consider  $\left\{i_{1}, \ldots, i_{s}\right\} \subset\{0,1, \ldots, n\}$ the set of ordered set of $s$ indices, that is, $i_{1}<i_{2}<\cdots<i_{s}$. Let us denote by I its power set (consisting of all of the $2^{s}$ subsets of the set), and by $J$ the set of all proper subsets. Given a $(2 n+2)$-tuple $(\mathbf{u}, \mathbf{v})=\left(u_{0}, \ldots, u_{n}, v_{0}, \ldots, v_{n}\right)$ of integers, let us  define inductively a set of $2^{s}$ positive integers, one for each ordered element of $I$, as follows:
$$
c_{\emptyset}=\operatorname{gcd}\left(u_{0}, \ldots, u_{n}\right),
$$
and if $\left\{i_{1}, \ldots, i_{s}\right\} \subset\{0,1, \ldots, n\}$ is ordered, then
$$
c_{i_{1}, \ldots, i_{s}}=\frac{\operatorname{gcd}\left(u_{0}, \ldots, \hat{u}_{i_{1}}, \ldots, \hat{u}_{i_{s}}, \ldots, u_{n}\right)}{\prod_{J} c_{j_{1}, \ldots j_{t}}} .
$$
Similarly, we also define a set of $2^{s}$ real numbers by
$$
k_{\emptyset}=\epsilon_{n+1},
$$
and
$$
k_{i_{1}, \ldots, i_{s}}=\epsilon_{n-s+1} \sum_{I}(-1)^{s-t} \frac{u_{j_{1}} \cdots u_{j_{t}}}{v_{j_{1}} \cdots v_{j_{t}} \operatorname{lcm}\left(u_{j_{1}}, \ldots, u_{j_{t}}\right)}
$$
where
$$
\epsilon_{n-s+1}= \begin{cases}0 & \text { if } n-s+1 \text { is even } \\ 1 & \text { if } n-s+1 \text { is odd }\end{cases}
$$
respectively. Finally, for any $j$ such that $1 \leq j \leq r=\left\lfloor\max \left\{k_{i_{1}, \ldots, i_{s}}\right\}\right\rfloor$, where $\lfloor x\rfloor$ is the greatest integer less than or equal to $x$, we set
$
d_{j}=\prod_{k_{i_{1}, \ldots, i_{s}} \geq j} c_{i_{1}, \ldots, i_{s}}.$
Then 

\begin{equation}
H_{n-1}\left(L_{f}, \mathbb{Z}\right)_{\text {tor }}=\mathbb{Z} / d_{1} \oplus \cdots \oplus \mathbb{Z} / d_{r} .\end{equation}

\end{conjecture}
\medskip

In  \cite{HM}, Hertling and Mase   showed that this conjecture is true  for the following cases:

\begin{enumerate}

\item {\it Chain type singularity:}  a quasihomogeneous singularity of the form
$$
f=f\left(x_{1}, \ldots, x_{n}\right)=x_{1}^{a_{1}+1}+\sum_{i=2}^{n} x_{i-1} x_{i}^{a_{i}}
$$
for some $n \in \mathbb{N}$ and some $a_{1}, \ldots, a_{n} \in \mathbb{N}$. 
\item {\it Cycle or loop type singularity:} a quasihomogeneous singularity of the form
$$
f=f\left(x_{1}, \ldots, x_{n}\right)=\sum_{i=1}^{n-1} x_{i}^{a_{i}} x_{i+1}+x_{n}^{a_{n}} x_{1}
$$
for some $n \in \mathbb{Z}_{\geq 2}$ and some $a_{1}, \ldots, a_{n} \in \mathbb{N}$ which satisfy for even $n$ neither $a_{j}=1$ for all even $j$ nor $a_{j}=1$ for all odd $j.$
\item {\it Thom-Sebastiani iterated sums of singularities of chain type or cyclic type.}  Recall, for  $f$ and $g$  singularities, the Thom-Sebastiani sum is given by  $$f+g=f\left(x_{1}, \ldots, x_{n_{f}}\right)+g\left(x_{n_{f}+1}, \ldots, x_{n_{f}+n_{g}}\right).$$ 
Notice that  {\it Brieskorn-Pham singularities}, or {\it BP singularities} $$f=f\left(x_{1}, \ldots, x_{n}\right)=\sum_{i=1}^{n} x_{i}^{a_{i}}$$ for some $n \in \mathbb{N}$ and some $a_{1}, \ldots, a_{n} \in \mathbb{Z}_{\geq 2}$ are special cases of chain type singularities. 
\end{enumerate}

Any iterated Thom-Sebastiani sum of chain type singularities and cycle type singularities are also called {\it invertible polynomials}.
In general, they can be written as $$f\left(x_1, \ldots, x_n\right)=\sum_{j=1}^n \prod_{i=1}^n x_i^{A_{i j}}$$ with $\operatorname{det}\left(A_{i j}\right) \neq 0$, then the Berglund-Hübsch transpose (BH-transpose for short)  \cite{BH}, is given by 
$$
f^T\left(x_1, \ldots, x_n\right)=\sum_{j=1}^n \prod_{i=1}^n x_i^{A_{j i}}, 
$$
that is,  $f^T\left(x_1, \ldots, x_n\right)$ is defined by transposing the exponent matrix $A=\left(A_{i j}\right)_{i, j}$ of the original polynomial.  It is well-known that the BH-transpose of an invertible polynomial remains invertible. 

In the early days of mirror symmetry, the BH-transpose rule was used to exhibit many examples of mirror pairs as hypersurfaces in weighted projective space.  Recall that mirror symmetry is a conjectured symmetry among Calabi-Yau varieties. In its classical formulation, it asserts that if $X$ is a Calabi-Yau variety of dimension $n$ then there exists another Calabi-Yau variety  $X^{\wedge}$ such that $$H^{p,q}(X,\mathbb{C})\cong H^{n-p,q}(X^{\wedge},\mathbb{C}).$$ 
The Berglund-H\"ubsch mirror symmetry construction takes $X_f$ a Calabi-Yau hypersurface in $\mathbb{P}(\mathbf{w})$ defined by the invertible polynomial $f$ of degree $d$ and form the exponent matrix $A_f.$       
Applying the BH-transpose rule one obtains the invertible polynomial $f^{T}$ with weight vector $\mathbf{w}_{T}$ and degree $d_T$ such that  $A^{T}\mathbf{w}_T=\mathbf{d}_T.$ This new quasi-homogeneous polynomial $f^{T}$ of degree $d_T$  defines the orbifold $X_{f^T}$ in  $\mathbb{P}(\mathbf{w}_T)$ which is again a Calabi-Yau hypersurface, where $A^{T}_{f}=A_{f^{T}}.$ Let us  represent this preliminary idea by the diagram:
 \begin{center}
\begin{tikzcd}
f \arrow[r] \arrow[d] & f^{T}\\
A_{f} \arrow[r,"T"] &A_{f}^{T}\arrow[u].\\
\end{tikzcd}
\end{center}

According to Berglund and H\"ubsch in \cite{BH}, the mirror to $X_f$ should be the quotient of the hypersurface $X_{f^T}$ by a subgroup $G$ of the diagonal symmetries of the hypersurface determined by $f$. Later, Krawitz generalized Berglund-H\"ubsch mirror symmetry construction introducing a {\it dual} group $G^T$ of $G$ which this time contains the diagonal symmetries induced by the ${\mathbb C}^*$-action on the weighted projective space (see details in \cite{Kr}).

 The important thing to note for our paper is that the BH-transpose rule preserves not only the Calabi-Yau condition but the Fano condition as well. Furthermore, as we will see, in many cases the rule also preserves the Sasaki-Einstein condition of the corresponding link,  provided this is a rational homology 7-sphere.

 \begin{lemma}[\cite{Go}]
 Let $f$ be an invertible polynomial defining a hypersurface $X_f$ in weighted projective space $\mathbb{P}(\mathbf{w}),$ let $X_{f^T}$ be the associated hypersurface in $\mathbb{P}\left(\mathbf{w}_T\right)$ defined by $f^T$ of degree $d_T$. Then if $X_f$ is a Fano hypersurface, then so is $X_{f^T}$ (to be more precise it is log-Fano).  Thus, if the link $L_f$ associated to $f$ admits positive Ricci curvature Sasaki metric then the link $L_{f^T}$ corresponding to its transpose also admits a Sasaki metric of positive Ricci curvature.
 \end{lemma}

It is natural to try to find conditions to extend this result at the level of  links that admit Sasaki-Einstein metrics. In order to achieve this, we focus on links that are rational homology 7-spheres. So far, through calculations given in \cite{BGN2} and  \cite{CL}  all the polynomials that produce rational homology 7-spheres and  admit Sasaki-Einstein metrics are polynomials of cycle type, chain type, or iterated Thom-Sebastiani sums of these types. Furthermore, all these examples satisfy one of the following conditions (proofs  of the following statements can be found in \cite{BGN2} or in the proof of Theorem 3.2 in \cite{CL}):

\begin{enumerate}
\item[A)] The link $L(\textbf{w}, d)$ is such that $\operatorname{gcd}(d,w_{i})=1$ for all $i=0,\dots 4.$ It follows that $$\mu + 1 = d(b_{n-1}+1)\,\, \mbox{and}\,\, H_{n-1}(L_f,\ZZ)_{tor}=\ZZ_{d}.$$ Thus,  for hypersurface singularities with  $\operatorname{gcd}\left(d, w_{i}\right)=1$ for all $i$ such that their corresponding links are rational homology spheres,  we have   $\mu=d-1.$ 
\item[B)] The link $L(\textbf{w}, d)$ is such that its weight vector $\textbf{w}=(w_{0},w_{1},w_{2},w_{3},w_{4})$ satisfies  $\textbf{w}=(w_{0},w_{1},w_{2},w_{3},w_{4})=(m_{2}v_{0},m_{2}v_{1},m_{2}v_{2},m_{3}v_{3},m_{3}v_{4})$ 
where the $v_i's$ are given as in Equation (1.0.2),   $\gcd(m_{2},m_{3})=1$ and $m_{2}m_{3}=d$ where  $m_2$ is odd and $m_3$ is even.  
One  obtains  the equality 
$$\operatorname{div}\Delta_f=\alpha({\bf w})\beta({\bf w})\Lambda_d+\beta({\bf w})\Lambda_{m_3}-\alpha({\bf w})\Lambda_{m_2}-\Lambda_1,$$ 
with the two positive  integers $\alpha(\textbf{w})$ and $\beta(\textbf{w})$ depending on the weights:
\begin{equation}
\alpha(\textbf{w})=\frac{m_{2}}{v_{3}v_{4}}-\frac{1}{v_{3}}-\frac{1}{v_{4}}
\end{equation}
and 

\begin{equation}
\beta(\textbf{w})=\left (\frac{m_3}{v_0v_1}-\frac{1}{v_1}-\frac{1}{v_0}\right )\left (\frac{m_3}{v_2}- 1\right ) + \frac{1}{v_2}.
\end{equation}
Furthermore, if  the link is a rational homology sphere,  then  $\beta({\bf w})=1$ and $H_{3}(L_f,\mathbb{Z})=(\mathbb{Z}_{m_{3}})^{\alpha(\textbf{w})+1}.$

\end{enumerate}


\section{Invariance of rational homology 7-spheres under BH-transpose rule}

Links of a weighted homogeneous singularity can be viewed as Seifert bundle $\pi:L\rightarrow X$ and these were studied by  Orlik and Wagreich in \cite{OW}. From the spectral sequence   
$$H^i\left(Z, R^j \pi_* \mathbb{Q}_L\right) \Rightarrow H^{i+j}(L, \mathbb{Q})$$ one obtains the following (see \cite{Kol1} for a neat argument): 
\begin{enumerate}
\item 
$\operatorname{dim} H^i(L, \mathbb{Q}) =\operatorname{dim} H^i(X, \mathbb{Q})-\operatorname{dim} H^{i-2}(X, \mathbb{Q}) \quad\hbox{for}\quad   i \leq \operatorname{dim} X,$ 
\item $\operatorname{dim} H^{i+1}(L, \mathbb{Q})  =\operatorname{dim} H^i(X, \mathbb{Q})-\operatorname{dim} H^{i+2}(X, \mathbb{Q})  \quad\hbox{for}\quad  i \geq \operatorname{dim} X.$ 
\end{enumerate}
Thus, $L$ is a rational homology sphere if and only if $X$ is a rational homology complex projective space. In the context of Sasaki-Einstein metrics, it is natural to ask whether the corresponding orbifolds, which must be Fano K\"ahler-Einstein, are birational to complex projective spaces or whether there are  toric actions associated with them that can be used to understand the special metric structures found on the total spaces of the fibration. Of interest to us in this article, are rational homology 7-spheres admitting Sasaki-Einstein metrics that arise as \emph{twins} \cite{BGN2}.

\begin{definition}
Let $L_{f}$ and $L_g$ be links of an isolated hypersurface singularity which are rational homology $(2n+1)$-spheres. The links $L_f$ and $L_g$ are said to be twins if they satisfy   $\mu(L_f)=\mu(L_g)$, $d_g=d_f$, and $H_{n}(L_f,\mathbb{Z})=H_{n}(L_g,\mathbb{Z}).$
\end{definition}
As mentioned in \cite{BGN2}, there is some evidence to conjecture that twins are homeomorphic or even diffeomorphic, however, as far the authors know, there is no proof for this claim.    
 The first examples of Sasaki-Einstein twins were given in \cite{BGN2} and more were constructed in \cite{CL}. Because Sasaki-Einstein links are contingent upon specific hypersurfaces in weighted projective space, the construction lends itself to techniques from an entirely different construction in the context of mirror symmetry: using ideas from the 
 BH-transpose rule from Berglund-H\"ubsch-Krawitz mirror symmetry, we are able to give a new framework for how the Sasaki-Einstein twins arise.

We have the following theorem, which includes as hypotheses certain arithmetic conditions that were present in the examples of rational homology 7-spheres found in \cite{BGN2} and \cite{CL}.

\begin{theorem} Consider the following data $(d, \mathbf{w}=(w_{0},w_{1},w_{2},w_{3},w_{4}) ),$ such that $d=w_{0}+w_{1}+w_{2}+w_{3}+w_{4}-1.$ Then 
\begin{enumerate}
\item Given the invertible polynomial of cycle type 
$$f=z_{4}z_{0}^{a_{0}}+z_{0}z_{1}^{a_{1}}+z_{1}z_{2}^{a_{2}}+z_{2}z_{3}^{a_{3}}+z_{3}z_{4}^{a_{4}}$$
of degree $d$ in the well-formed projective space $\mathbb{P}(\mathbf{w})$  such that $\gcd(d,w_{i})=1.$ Suppose its corresponding link $L_f$ is a rational homology sphere.
Then the link $L_{f_{T}}$ corresponding to the  BH-transpose rule of $f$  is a rational homology sphere.  
\item Given any of the following invertible polynomials of degree $d$ in the the well-formed  projective space $\mathbb{P}(\mathbf{w})$
\begin{align*}
\text{chain-cycle: }f &=z_{0}^{a_{0}}+z_{0}z_{1}^{a_{1}}+z_{4}z_{2}^{a_{2}}+z_{2}z_{3}^{a_{3}}+z_{3}z_{4}^{a_{4}} \\
\text{BP-ycle: } f &=z_0^{a_0}+z_1^{a_1}+z_4 z_2^{a_2}+z_2 z_3^{a_3}+z_3 z_4^{a_4}\\
\text{cycle-cycle: } f &=z_{1}z_{0}^{a_{0}}+z_{0}z_{1}^{a_{1}}+z_{4}z_{2}^{a_{2}}+z_{2}z_{3}^{a_{3}}+z_{3}z_{4}^{a_{4}}
\end{align*}
such that  
$\mathbf{w}=(w_{0},w_{1},w_{2},w_{3},w_{4})=(m_{3}v_{0},m_{3}v_{1},m_{2}v_{2},m_{2}v_{3},m_{2}v_{4}),$
with $\gcd(m_{2},m_{3})=1$ and $d=m_3m_2.$ If  the corresponding link $L_f$ is a rational homology sphere, then the link $L_{f_{T}}$ associated to  the  BH-transpose rule of $f$  is a rational homology sphere. 
\end{enumerate}

\end{theorem}
\noindent{\it Proof.} The proof follows from Cases I, II, III and IV presented below. Cases I, II and III showed in particular that the BH-transpose rule only produces twins. Case IV is more involved and interesting, the BH-transpose rule does not preserve neither torsion nor Milnor number.  

\hfill$\square$
\medskip

\subsection{Case I: cycle type}
Given the invertible polynomial of cycle type 
$$f=z_{4}z_{0}^{a_{0}}+z_{0}z_{1}^{a_{1}}+z_{1}z_{2}^{a_{2}}+z_{2}z_{3}^{a_{3}}+z_{3}z_{4}^{a_{4}}$$
of degree $d$ in the projective space $\mathbb{P}(\textbf{w})$, where  $\textbf{w}=(w_{0},w_{1},w_{2},w_{3},w_{4})$ such that $d=w_{0}+w_{1}+w_{2}+w_{3}+w_{4}-1$. The associated  matrix of exponentials for $f$ is given by 
\begin{equation*}
    A=\begin{bmatrix}
        a_{0} & 0 & 0 & 0 & 1 \\
        1 & a_{1} & 0 & 0 & 0 \\
        0 & 1 & a_{2} & 0 & 0 \\
        0 & 0 & 1 & a_{3} & 0 \\
        0 & 0 & 0 & 1 & a_{4}
    \end{bmatrix}
\end{equation*}
Here, using the BH-transpose rule,  we associate the invertible polynomial 
$$f_{T} = z_{0}^{a_{0}}z_{1}+ z_{1}^{a_{1}}z_{2}+z_{2}^{a_{2}}z_{3}+z_{3}^{a_{3}}z_{4}+z_{4}^{a_{4}}z_{0} .$$ 
The aim is to find  weights $\textbf{w}_{T}=(\tilde{w}_{0},\tilde{w}_{1},\tilde{w}_{2},\tilde{w}_{3},\tilde{w}_{4})$ that satisfy the matrix equation
\begin{equation}
A^{T}\Tilde{W}=D
\end{equation}
where $\Tilde{W}=(\Tilde{w}_{0},\Tilde{w}_{1},\Tilde{w}_{2},\Tilde{w}_{3},\Tilde{w}_{4})^{T}$ and $D=(d,d,d,d,d)^{T}$.

\noindent{\bf First step:} We will prove that the column matrix  $W=(\tilde{w}_{0},\tilde{w}_{1},\tilde{w}_{2},\tilde{w}_{3},\tilde{w}_{4})^{T}$obtained as solution of (3.1.1) has as entries integer numbers and the index of $f_T$ equals 1. Indeed, for the polynomial $f$ we have that the the matrix equation:
\begin{equation}
    AW=D 
\end{equation}
where $W=(w_{0},w_{1},w_{2},w_{3},w_{4})^{T}$  and $D=(d,d,d,d,d)^{T}$ has as solution $A.$ Thus, we have  
\begin{equation}a_{0}w_{0}+w_{4}=d, \ \ a_{1}w_{1}+w_{0}= d, \ \ a_{2}w_{2}+w_{1}=d, \ \ a_{3}w_{3}+w_{2}=d, \ \ a_{4}w_{4}+w_{3}=d
\end{equation}
Solving these equations for $w_{4}$, we obtain
\begin{equation}
w_{4}=\dfrac{d}{1+a_{0}a_{1}a_{2}a_{3}a_{4}}(1-a_{0}+a_{0}a_{1}-a_{0}a_{1}a_{2}+a_{0}a_{1}a_{2}a_{3}).
\end{equation}
Since the corresponding link for $f$ is a rational homology sphere and $\gcd(d,w_{i})=1$  we have $\mu=d-1.$   Furthermore, by (3.1.3) we have
$$a_{0}a_{1}a_{2}a_{3}a_{4}=\left(\dfrac{d-w_{4}}{w_{0}}\right)\left(\dfrac{d-w_{0}}{w_{1}}\right)\left(\dfrac{d-w_{1}}{w_{2}}\right)\left(\dfrac{d-w_{2}}{w_{3}}\right)\left(\dfrac{d-w_{3}}{w_{4}}\right) =\mu.$$
Replacing in (3.1.4), we obtain the
following equation for the weights, using subscripts mod 5:
\begin{equation}
    w_i=1-a_{i+1}+a_{i+1}a_{i+2}-a_{i+1}a_{i+2}a_{i+3}+a_{i+1}a_{i+2}a_{i+3}a_{i+4}.
\end{equation}
Repeating these steps to solve the matrix equation (3.1.1), we obtain the equation using subscripts mod 5:

\begin{equation}
\tilde{w}_{i}=1-a_{i+4}+a_{i+4}a_{i+3}-a_{i+4}a_{i+3}a_{i+2}+a_{i+4}a_{i+3}a_{i+2}a_{1+1}.
\end{equation}

Thus $\tilde{w}_{i}\in\mathbb{Z}$.  Furthermore, from (3.1.5) and (3.1.6), we notice that the weights verify,
\begin{equation}
\tilde{w}_{0}+\tilde{w}_{1}+\tilde{w}_{2}+\tilde{w}_{3}+\tilde{w}_{4}=w_{0}+w_{1}+w_{2}+w_{3}+w_{4} = d+1.
\end{equation}
 
 \noindent{\bf Second step:}  We will show that the third Betti number and the torsion of $H_3$ do not change under BH-transpose rule. To show this, first, let us see that $\gcd(d,\tilde{w}_{i})=1$, for all $i$. Indeed, from (2.1.1) we have the equation, using indices mod 5:
\begin{equation}
\tilde{w}_{i}=d-a_{i+4}\tilde{w}_{i+4}.
\end{equation}

If we suppose that $\gcd(d,\tilde{w}_{i_{0}})=l\neq1$, for some $i_{0}$, then by (3.1.8) we have $\gcd(d,\tilde{w}_{i})=l$, for all $i$. Indeed,  let us put $\gcd(d, \tilde{w}_{i})=\alpha_i.$ From Equation (3.1.8) we have  $\gcd(d, \tilde{w}_{0})=\gcd(d, d-a_4 \tilde{w}_{0})=\gcd(d, a_4\tilde{w}_{4}).$ We conclude that  $\alpha_4/ \alpha_0.$
 Using the same argument as above, if  $\alpha_1=\gcd(d, \tilde{w}_{1})=\gcd(d, d-a_o \tilde{w}_{0})=\gcd(d, a_0\tilde{w}_{0}),$ thus $\alpha_0/ \alpha_1.$ Repeating this argument, we conclude that $\alpha_0/ \alpha_1,\alpha_1/ \alpha_2,\alpha_2/ \alpha_3,\alpha_3/ \alpha_4$ and $\alpha_4/ \alpha_0.$ It follows  
 that $\gcd(d,\tilde{w}_{i})=\alpha_0=l$, for all $i$ as claimed. Moreover, from Equation (3.1.7) we have that $\gcd(d, \tilde{w}_{0}+\tilde{w}_{1}+\tilde{w}_{2}+\tilde{w}_{3}+\tilde{w}_{4})=\gcd(d, d+1)=1.$ It follows that $l=1.$

On the other hand, from (3.1.8) we notice that 
\begin{equation}
    a_{0}a_{1}a_{2}a_{3}a_{4}=\left(\dfrac{d-\tilde{w}_{1}}{\tilde{w}_{0}}\right)\left(\dfrac{d-\tilde{w}_{2}}{\tilde{w}_{1}}\right)\left(\dfrac{d-\tilde{w}_{3}}{\tilde{w}_{2}}\right)\left(\dfrac{d-\tilde{w}_{4}}{\tilde{w}_{3}}\right)\left(\dfrac{d-\tilde{w}_{0}}{\tilde{w}_{4}}\right) = \tilde{\mu}.
\end{equation}
Since $\mu=a_{0}a_{1}a_{2}a_{3}a_{4}$ we have  $\mu=\tilde{\mu}$. Using this result and the argument given in A), we have that 
$\tilde{b}_{3}=b_{3}=0$ and the  torsion remains the same.

\subsection{Case II: Brieskorn-Pham-cycle}

Even though, this is a particular case of the Chain-Cycle case, we prefer to work the details separately.

We will work with hypersurfaces $X_{f}$ where $f$ can be written as
$$\mbox{BP-cycle: } f=z_{0}^{a_{0}}+z_{1}^{a_{1}}+z_{4}z_{2}^{a_{2}}+z_{2}z_{3}^{a_{3}}+z_{3}z_{4}^{a_{4}}.$$
 Then we can associate its exponent matrix 
\begin{equation*}
    A =\begin{bmatrix}
      a_{0} & 0 & 0 & 0 & 0 \\
      0 & a_{1} & 0 & 0 & 0 \\
      0 & 0 & a_{2} & 0 & 1 \\
      0 & 0 & 1 & a_{3} & 0 \\
      0 &  0 & 0 & 1 & a_{4} 
    \end{bmatrix}
\end{equation*}
Using the  BH transpose rule  
we obtain the polynomial $\Tilde{f}_{T}=z_{0}^{a_{0}}+z_{1}^{a_{1}}+z_{3}z_{2}^{a_{2}}+z_{4}z_{3}^{a_{3}}+z_{2}z_{4}^{a_{4}}$ with weights $\textbf{w}_{T}=(\Tilde{w}_{0},\Tilde{w}_{1},\Tilde{w}_{2},\Tilde{w}_{3},\Tilde{w}_{4})$. We will calculate the vector weight $\textbf{w}_{T}$ solving the matrix equation
\begin{equation}
A^{T}\Tilde{W}=D
\end{equation}
where $\Tilde{W}=(\Tilde{w}_{0},\Tilde{w}_{1},\Tilde{w}_{2},\Tilde{w}_{3},\Tilde{w}_{4})^{T}$ and $D=(d,d,d,d,d)^{T}$.

\noindent\textbf{First step:} Let us show that the entries of the  column matrix obtained as solution of (3.2.1) are integers and the index of $f^T$ is 1. 

The weight vector $\textbf{w}=(w_{0},w_{1},w_{2},w_{3},w_{4})=(m_{3}v_{0},m_{3}v_{1},m_{2}v_{2},m_{2}v_{3},m_{2}v_{4})$ satisfies 
   \begin{equation}
    AW=D 
\end{equation}
where $W=(w_{0},w_{1},w_{2},w_{3},w_{4})^{T}$  and $D=(d,d,d,d,d)^{T}$ has as solution $A.$
We obtain
   
   \begin{align}
        w_{0} &=\dfrac{1}{a_{0}}d, w_{1}=\dfrac{1}{a_{1}}d, \nonumber\\
        w_{2} &=\dfrac{d}{a_{2}a_{3}a_{4}+1}(a_{4}a_{3}-a_{3}+1),\nonumber\\ 
        w_{3} &= \dfrac{d}{a_{2}a_{3}a_{4}+1}(a_{2}a_{4}-a_{4}+1),\\ 
        w_{4} &= \dfrac{d}{a_{2}a_{3}a_{4}+1}(a_{2}a_{3}-a_{2}+1). \nonumber
    \end{align}
    Now we compute $a_{2}a_{3}a_{4}$. From (3.2.1) we have 
    \begin{equation}
        a_{2}w_{2}+w_{4}=d, \ \ a_{3}w_{3}+w_{2}=d \ \mbox{ and  } \ a_{4}w_{4}+w_{3}=d.
    \end{equation}
       Then, from the assumptions on the weight vector ${\bf w}$ we have 
       \begin{align*}
           a_{2}a_{3}a_{4} & = \left(\dfrac{d-w_{4}}{w_{2}}\right)\left(\dfrac{d-w_{2}}{w_{3}}\right)\left(\dfrac{d-w_{3}}{w_{4}}\right) \\
           & = \left(\dfrac{m_{3}-v_{4}}{v_{2}}\right)\left(\dfrac{m_{3}-v_{2}}{v_{3}}\right)\left(\dfrac{m_{3}-v_{3}}{v_{4}}\right) \\
           & = \dfrac{m_{3}^{3}-(v_{2}+v_{3}+v_{4})m_{3}^{2}+(v_{2}v_{3}+v_{2}v_{4}+v_{3}v_{4})m_{3}-v_{2}v_{3}v_{4}}{v_{2}v_{3}v_{4}} \\
           & = m_{3}\left(\dfrac{m_{3}^{2}-(v_{2}+v_{3}+v_{4})m_{3}+v_{2}v_{3}+v_{2}v_{4}+v_{3}v_{4}}{v_{2}v_{3}v_{4}}\right)-1.
       \end{align*}
      Since the corresponding link is a rational homology sphere, it follows  from Equation (2.2.8) that 
    $$\beta(\textbf{w})=\dfrac{m_{3}^{2}-(v_{2}+v_{3}+v_{4})m_{3}+v_{2}v_{3}+v_{2}v_{4}+v_{3}v_{4}}{v_{2}v_{3}v_{4}}=1.$$
    Substituting this equality in the previous equation,  we obtain
    \begin{equation}
        a_{2}a_{3}a_{4}+1=m_{3}.
    \end{equation}
    Replacing (3.2.5) in (3.2.3), we have
    \begin{equation}
        w_{2}=m_{2}(a_{4}a_{3}-a_{3}+1), \quad   w_{3} = m_{2}(a_{2}a_{4}-a_{4}+1) \ \ \mbox{and} \ \ w_{4}=m_{2}(a_{2}a_{3}-a_{2}+1).    
    \end{equation}
        
    On the other hand, solving (3.2.1) and using (3.2.6) we obtain the following expressions for $\tilde{w}_{i}:$ 

\begin{equation}
   \begin{split}
   \tilde{w}_{0} & =\dfrac{1}{a_{0}}d=w_{0}, \tilde{w}_{1}=\dfrac{1}{a_{1}}d=w_{1},  \\
   \tilde{w}_{2}&=m_{2}(a_{4}a_{3}-a_{4}+1),\\
   \tilde{w}_{3} &= m_{2}(a_{2}a_{4}-a_{2}+1),\\ 
   \tilde{w}_{4} & = m_{2}(a_{2}a_{3}-a_{3}+1).
   \end{split}
\end{equation}
Thus $\tilde{w}_{i}\in\mathbb{Z}$. Furthermore, from (3.2.6) and (3.2.7) we obtain $$\tilde{w}_{0}+\tilde{w}_{1}+\tilde{w}_{2}+\tilde{w}_{3}+\tilde{w}_{4}-1=w_{0}+w_{1}+w_{2}+w_{3}+w_{4}-1=d.$$

\noindent\textbf{Second step:} We will show that $\tilde{b}_{3}=0$ for $\textbf{w}_{T}=(\Tilde{w}_{0},\Tilde{w}_{1},\Tilde{w}_{2},\Tilde{w}_{3},\Tilde{w}_{4})$.

By Equation (3.2.7), we can express
    $$\textbf{w}_{T}=(\Tilde{w}_{0},\Tilde{w}_{1},\Tilde{w}_{2},\Tilde{w}_{3},\Tilde{w}_{4})=(m_{3}\Tilde{v}_{0},m_{3}\Tilde{v}_{1},m_{2}\Tilde{v}_{2},m_{2}\Tilde{v}_{3},m_{2}\Tilde{v}_{4}),$$
    where $\tilde{v}_{0}=v_{0}, \tilde{v}_{1}=v_{1}, \tilde{v}_{2}=a_{4}a_{3}-a_{4}+1, \tilde{v}_{3}=a_{2}a_{4}-a_{2}+1 $ and $\tilde{v}_{4}=a_{2}a_{3}-a_{3}+1$.

   From (3.2.4), we have
   \begin{align}
       a_{2}w_{2}+w_{4}=d & \Rightarrow a_{2}m_{2}v_{2}+m_{2}v_{4} =m_{2}m_{3} \Rightarrow a_{2}=\dfrac{m_{3}-v_{4}}{v_{2}}, \\ 
       a_{3}w_{3}+w_{2}=d & \Rightarrow a_{3}m_{2}v_{3}+m_{2}v_{2} =m_{2}m_{3} \Rightarrow a_{3}=\dfrac{m_{3}-v_{2}}{v_{3}}, \\ 
       a_{4}w_{4}+w_{3}=d & \Rightarrow a_{4}m_{2}v_{4}+m_{2}v_{3} =m_{2}m_{3} \Rightarrow a_{4}=\dfrac{m_{3}-v_{3}}{v_{4}}, 
   \end{align}

 Similarly, from (3.2.1) we obtain

   \begin{align}
       a_{2}\tilde{w}_{2}+\tilde{w}_{3}=d & \Rightarrow a_{2}\tilde{v}_{2}+\tilde{v}_{3} =m_{3}, \\ 
       a_{3}\tilde{w}_{3}+\tilde{w}_{4}=d & \Rightarrow a_{3}\tilde{v}_{3}+\tilde{v}_{4} =m_{3}, \\ 
       a_{4}\tilde{w}_{4}+\tilde{w}_{2}=d & \Rightarrow a_{4}\tilde{v}_{4}+\tilde{v}_{2} =m_{3}.  
   \end{align}
   Replacing equations (3.2.8)-(3.2.10) in equations (3.2.11)-(3.2.13) we have 
   $$
\left\{
\begin{array}{rcl}
\left(\dfrac{m_{3}-v_{4}}{v_{2}}\right)\tilde{v}_{2}=m_{3}-\tilde{v}_{3}, \\
\left(\dfrac{m_{3}-v_{2}}{v_{3}}\right)\tilde{v}_{3}=m_{3}-\tilde{v}_{4}, \\
\left(\dfrac{m_{3}-v_{3}}{v_{4}}\right)\tilde{v}_{4}=m_{3}-\tilde{v}_{2}.
\end{array}
\right. 
$$
Multiplying the expressions above, we have the following equality
$$\dfrac{(m_{3}-v_{2})(m_{3}-v_{3})(m_{3}-v_{4})}{v_{2}v_{3}v_{4}}=\dfrac{(m_{3}-\tilde{v}_{2})(m_{3}-\tilde{v}_{3})(m_{3}-\tilde{v}_{4})}{\tilde{v}_{2}\tilde{v}_{3}\tilde{v}_{4}}.$$
Simplifying, we have
{\small{
\begin{equation}
    \dfrac{m_{3}^{2}-(v_{2}+v_{3}+v_{4})m_{3}+v_{2}v_{3}+v_{2}v_{4}+v_{3}v_{4}}{v_{2}v_{3}v_{4}}=\dfrac{m_{3}^{2}-(\tilde{v}_{2}+\tilde{v}_{3}+\tilde{v}_{4})m_{3}+\tilde{v}_{2}\tilde{v}_{3}+\tilde{v}_{2}\tilde{v}_{4}+\tilde{v}_{3}\tilde{v}_{4}}{\tilde{v}_{2}\tilde{v}
    _{3}\tilde{v}_{4}}
\end{equation}
}}
Using again Equation (2.2.8), we have 
$$\dfrac{m_{3}^{2}-(v_{2}+v_{3}+v_{4})m_{3}+v_{2}v_{3}+v_{2}v_{4}+v_{3}v_{4}}{v_{2}v_{3}v_{4}}=\beta(\textbf{w})=1,$$ since $L_f$ is a rational homology sphere.
Thus $\beta(\tilde{\textbf{w}})=1$ and $\tilde{b}_{3}=0$.
\medskip

 \noindent\textbf{Third step:} We will prove that  the torsion remains invariant under the BH-transpose rule. 

 \medskip

 For the weights $\textbf{w}=(w_{0},w_{1},w_{2},w_{3},w_{4})=(m_{3}v_{0},m_{3}v_{1},m_{2}v_{2},m_{2}v_{3},m_{2}v_{4})$, we define $$u_{i}=\dfrac{d}{\gcd(d,w_{i})}.$$ 
 Since $(w_{0},w_{1},w_{2},w_{3},w_{4})$ defines a well-formed weighted projective space, we obtain $u_{0}=u_{1}=m_{2}.$ Now, we will show that $u_{2}=u_{3}=u_{4}=m_{3}.$ Suppose that 
 $$\gcd(d,w_{i})=gcd(m_{2}m_{3},m_{2}v_{i})\neq m_{2}, \mbox{ for some }i=2,3,4.$$
 This implies that $\gcd(m_{3},v_{i})\neq 1$. Assume without loss of generality that $i=2$. Then using (3.2.8) we obtain that $\gcd(m_{3},v_{2},v_{4})\neq 1$. However, this is not possible since $(w_{0},w_{1},w_{2},w_{3},w_{4})$ defines a well-formed weighted projective space. Therefore $\gcd(d,w_{i})=m_{2}$ for all $i=2,3,4$. Now, using the formula for torsion given in Conjecture 2.2, we obtain $c_{\emptyset}=1, c_{i}=1,$ for all $i$. On the other hand, $c_{01}=m_{3}$ and $c_{ij}=1$  otherwise. Also we obtain $c_{234}=m_{2}$ and $c_{ijk}=1$ otherwise. Hence,  we just need calculate the values for $k_{01}$ and $k_{234}$. By formula, $k_{234}=0$ and 
 $$k_{01}=\dfrac{m_{2}}{v_{0}v_{1}}-\dfrac{1}{v_{0}}-\dfrac{1}{v_{1}}+1 = \alpha(\textbf{w})+1.$$
 Therefore, from the discussion  given in B), we have  $H_{3}(L_f,\mathbb{Z})=(\mathbb{Z}_{m_{3}})^{\alpha(\textbf{w})+1}.$ As $\tilde{v}_{0}=v_{0}$ and $\tilde{v}_{1}=v_{1}$, then $H_{3}(L_{f^T}, \mathbb{Z})=(\mathbb{Z}_{m_{3}})^{\alpha(\textbf{w})+1}$. 
\medskip

 \noindent\textbf{Fourth step:} We will prove that $L_{f}$ and $L_{f^T}$ have the same Milnor number.

 \medskip

From equations (3.2.5) to (3.2.10) we have

$$a_{2}a_{3}a_{4}= \left(\dfrac{d-w_{4}}{w_{2}}\right)\left(\dfrac{d-w_{2}}{w_{3}}\right)\left(\dfrac{d-w_{3}}{w_{4}}\right)=\left(\dfrac{d-\tilde{w}_{3}}{\tilde{w}_{2}}\right)\left(\dfrac{d-\tilde{w}_{4}}{\tilde{w}_{3}}\right)\left(\dfrac{d-\tilde{w}_{2}}{w_{4}}\right).$$
 Furthermore, as $w_{0}=\tilde{w}_{0}$ and $w_{1}=\tilde{w}_{1}$, then  $\mu(M)=\mu(\tilde{M}).$
\medskip

\subsection{Case III: cycle-cycle}

The polynomials of  cycle-cycle type are of the form 

$$\mbox{Cycle - Cycle: } f=z_{1}z_{0}^{a_{0}}+z_{0}z_{1}^{a_{1}}+z_{4}z_{2}^{a_{2}}+z_{2}z_{3}^{a_{3}}+z_{3}z_{4}^{a_{4}}.$$



The result is obtained in a similar way to the previous case. Here, we take $\tilde{w}_{0}=w_{0}$, $\tilde{w}_{1}=w_{1}$ and $\tilde{w}_{2},\tilde{w}_{3}$ and $\tilde{w}_{4}$ as in (3.2.4).

\subsection{Case IV: chain-cycle} 

Consider the  hypersurfaces $X_{f}$ where $f$ can be written as
$$\mbox{chain-cycle: } f=z_{0}^{a_{0}}+z_{0}z_{1}^{a_{1}}+z_{4}z_{2}^{a_{2}}+z_{2}z_{3}^{a_{3}}+z_{3}z_{4}^{a_{4}}.$$
Then we can associate its  exponent matrix 
\begin{equation*}
    A =\begin{bmatrix}
      a_{0} & 0 & 0 & 0 & 0 \\
      1 & a_{1} & 0 & 0 & 0 \\
      0 & 0 & a_{2} & 0 & 1 \\
      0 & 0 & 1 & a_{3} & 0 \\
      0 &  0 & 0 & 1 & a_{4}
    \end{bmatrix}
\end{equation*}
Using the method of BH transpose rule, we have its corresponding  polynomial $\Tilde{f}_{T}=z_{0}^{a_{0}}z_{1}+z_{1}^{a_{1}}+z_{3}z_{2}^{a_{2}}+z_{4}z_{3}^{a_{3}}+z_{2}z_{4}^{a_{4}}$ with weights $\textbf{w}_{T}=(\Tilde{w}_{0},\Tilde{w}_{1},\Tilde{w}_{2},\Tilde{w}_{3},\Tilde{w}_{4})$. We will calculate the vector weight $\textbf{w}_{T}$ solving the following matrix equation:
\begin{equation}
A^{T}\Tilde{W}=\Tilde{D},
\end{equation}
where $\Tilde{W}=(\Tilde{w}_{0},\Tilde{w}_{1},\Tilde{w}_{2},\Tilde{w}_{3},\Tilde{w}_{4})^{T}$ and $\Tilde{D}=(\Tilde{d},\Tilde{d},\Tilde{d},\Tilde{d},\Tilde{d})^{T}$.

\vspace{0.2cm}

This system can be written as
 \begin{equation}
        \begin{bmatrix}
          a_{0} & 1 & 0 & 0 & 0 & -1 \\
      0 & a_{1} & 0 & 0 & 0 & -1\\
      0 & 0 & a_{2} & 1 & 0 & -1\\
      0 & 0 & 0 & a_{3} & 1 & -1\\
      0 &  0 & 1 & 0 & a_{4} & -1  
        \end{bmatrix} \begin{bmatrix}
            \Tilde{w}_{0} \\
            \Tilde{w}_{1} \\
            \Tilde{w}_{2} \\
            \Tilde{w}_{3} \\
            \Tilde{w}_{4} \\
            \Tilde{d}
        \end{bmatrix} = \begin{bmatrix}
            0 \\ 
            0 \\
            0 \\
            0 \\
            0
        \end{bmatrix}
    \end{equation}
Since $A^{T}$ is an invertible matrix, we have that the solution space is one dimensional. Notice that if we consider a nonzero vector solution with integer entries $(\Tilde{w}_{0}, \Tilde{w}_{1},\Tilde{w}_{2},\Tilde{w}_{3},\Tilde{w}_{4},\Tilde{d})$, then the vector $$(w_{0}',w_{1}',w_{2}',w_{3}',w_{4}',d')=\dfrac{1}{\gcd(\Tilde{w}_{0},\Tilde{w}_{1},\Tilde{w}_{2},\Tilde{w}_{3},\Tilde{w}_{4},\Tilde{d})}(\Tilde{w}_{0},\Tilde{w}_{1},\Tilde{w}_{2},\Tilde{w}_{3},\Tilde{w}_{4},\Tilde{d})$$
is also a solution vector with integer entries. Now we define the numbers 

$$\Tilde{u}_{i}=\dfrac{\Tilde{d}}{\gcd(\Tilde{d},\Tilde{w}_{i
})}, \ \ \  \Tilde{v}_{i}=\dfrac{\Tilde{w}_{i}}{\gcd(\Tilde{d},\Tilde{w}_{i})}, \ \ \ u_{i}'=\dfrac{d'}{\gcd(d',w_{i}')} \ \mbox{ and } \ v_{i}'=\dfrac{w_{i}'}{\gcd(d',w_{i}')}.$$
These numbers verify that $\Tilde{u}_{i}=u_{i}'$ and $\Tilde{v}_{i}=v_{i}'$ for all $i=0,\dots,4$.  Thus, it is enough to  exhibit a vector solution with integer entries for the system (3.4.2) since the $u_i's$ do not vary and  these numbers are  the ones let us use of Orlik's formula.

 We put $\Tilde{d}$ as:
\begin{equation}
\Tilde{d}=d(m_{2}-1)=m_{3}m_{2}(m_{2}-1).
\end{equation}
\noindent\textbf{First step:} We will show that the entries of the column matrix obtained as solution of (3.4.1) are integers.

Since  $\textbf{w}=(w_{0},w_{1},w_{2},w_{3},w_{4})=(m_{3}v_{0},m_{3}v_{1},m_{2}v_{2},m_{2}v_{3},m_{2}v_{4}),$ we have 
\begin{equation}
AW=D,
\end{equation}
 where $W=(w_{0},w_{1},w_{2},w_{3},w_{4})^{T}$  and $D=(d,d,d,d,d)^{T}$ has as solution $A.$

    Solving this, we obtain
    \begin{equation}
        w_{0}=\dfrac{d}{a_{0}}, w_{1}=\dfrac{d-w_{0}}{a_{1}}
    \end{equation}
   and 
\begin{align}
        w_{2}&=\dfrac{d}{a_{2}a_{3}a_{4}+1}(a_{4}a_{3}-a_{3}+1),\nonumber\\ 
        w_{3} &= \dfrac{d}{a_{2}a_{3}a_{4}+1}(a_{2}a_{4}-a_{4}+1),\\
         w_{4} &= \dfrac{d}{a_{2}a_{3}a_{4}+1}(a_{2}a_{3}-a_{2}+1).\nonumber
\end{align}
 We notice  that for chain-cycle polynomials, $v_{0}=\frac{w_0}{\gcd(w_0, d)}=1$ which implies that $w_0=m_3.$ Thus, from (3.4.5) we obtain 
    \begin{equation}
        a_{0}=\dfrac{d}{w_{0}}=m_{2} \mbox{ and } a_{1}=\dfrac{d-w_{0}}{w_{1}}=\dfrac{m_{3}m_{2}-m_{3}}{m_{3}v_{1}}=\dfrac{m_{2}-1}{v_{1}}.
    \end{equation}
    On the other hand, we compute $a_{2}a_{3}a_{4}$. From (3.4.4) we have 
    \begin{equation}
        a_{2}w_{2}+w_{4}=d, \ \ a_{3}w_{3}+w_{2}=d \ \mbox{ and  } \ a_{4}w_{4}+w_{3}=d.
    \end{equation}
       Then
       \begin{align*}
           a_{2}a_{3}a_{4} & = \left(\dfrac{d-w_{4}}{w_{2}}\right)\left(\dfrac{d-w_{2}}{w_{3}}\right)\left(\dfrac{d-w_{3}}{w_{4}}\right) \\
           & = \left(\dfrac{m_{3}-v_{4}}{v_{2}}\right)\left(\dfrac{m_{3}-v_{2}}{v_{3}}\right)\left(\dfrac{m_{3}-v_{3}}{v_{4}}\right) \\
           & = \dfrac{m_{3}^{3}-(v_{2}+v_{3}+v_{4})m_{3}^{2}+(v_{2}v_{3}+v_{2}v_{4}+v_{3}v_{4})m_{3}-v_{2}v_{3}v_{4}}{v_{2}v_{3}v_{4}} \\
           & = m_{3}\left(\dfrac{m_{3}^{2}-(v_{2}+v_{3}+v_{4})m_{3}+v_{2}v_{3}+v_{2}v_{4}+v_{3}v_{4}}{v_{2}v_{3}v_{4}}\right)-1. 
       \end{align*}
         
From Equation (2.2.8),  we know   $$\beta(\textbf{w})=\dfrac{m_{3}^{2}-(v_{2}+v_{3}+v_{4})m_{3}+v_{2}v_{3}+v_{2}v_{4}+v_{3}v_{4}}{v_{2}v_{3}v_{4}}=1.$$
    Thus 
    \begin{equation}
        a_{2}a_{3}a_{4}+1=m_{3}.
    \end{equation}
 
    Returning to the equation (3.4.1) we have the following system:

    $$
\left\{
\begin{array}{rcl}
a_{0}\Tilde{w_{0}}+\Tilde{w_{1}}=\Tilde{d} \\
a_{1}\Tilde{w_{1}}=\Tilde{d} \\
a_{2}\Tilde{w_{2}}+\Tilde{w_{3}}=\Tilde{d} \\
a_{3}\Tilde{w_{3}}+\Tilde{w_{4}}=\Tilde{d}\\
a_{4}\Tilde{w_{4}}+\Tilde{w_{2}}=\Tilde{d}.
\end{array}
\right. 
$$
Solving the above system and using (3.4.3), (3.4.7) and (3.4.9) we obtain

$$\Tilde{w_{1}}=\dfrac{\Tilde{d}}{a_{1}}=m_{3}m_{2}v_{1},$$
$$\Tilde{w_{0}}=\dfrac{\Tilde{d}-\Tilde{w_{1}}}{a_{0}}=\dfrac{m_{3}m_{2}(m_{2}-1)-m_{3}m_{2}v_{1}}{m_{2}}=m_{3}v_{1}(a_{1}-1),$$
$$\Tilde{w}_{2}=\dfrac{\Tilde{d}(1-a_{4}+a_{4}a_{3})}{1+a_{2}a_{3}a_{4}}=m_{2}(m_{2}-1)(1-a_{4}+a_{3}a_{4}),$$
$$\Tilde{w}_{3}=\dfrac{\Tilde{d}(1-a_{2}+a_{2}a_{4})}{1+a_{2}a_{3}a_{4}}=m_{2}(m_{2}-1)(1-a_{2}+a_{2}a_{4}),$$
$$\Tilde{w}_{2}=\dfrac{\Tilde{d}(1-a_{3}+a_{2}a_{3})}{1+a_{2}a_{3}a_{4}}=m_{2}(m_{2}-1)(1-a_{3}+a_{2}a_{3}).$$
It follows that  $\Tilde{w}_{0},\Tilde{w}_{1},\Tilde{w}_{2},\Tilde{w}_{3},\Tilde{w}_{4}\in\mathbb{Z}$.
\medskip

\noindent\textbf{Second step:} Next we will show that $b_{3}=0$ for $\textbf{w}_{T}=(\Tilde{w}_{0},\Tilde{w}_{1},\Tilde{w}_{2},\Tilde{w}_{3},\Tilde{w}_{4})$. Here we use the formula  (2.2.5) given by Orlik to calculate the $(n-1)^{s t}$ Betti number.

Let us compute each $\gcd(\Tilde{d},\Tilde{w}_{i}).$ First, for $i=0$ we have
\begin{align*}
\gcd(\Tilde{d},\Tilde{w}_{0})&=\gcd(m_{3}m_{2}(m_{2}-1),m_{3}v_{1}(a_{1}-1))\\
 &=m_{3}v_{1}\gcd(m_{2}a_{1},a_{1}-1)\\
& =m_{3}v_{1}\gcd(m_{2},a_{1}-1).
\end{align*}

We write $a_{1}-1$ as
$$a_{1}-1=\dfrac{m_{2}-1}{v_{1}}-1=\dfrac{m_{2}-1-v_{1}}{v_{1}}.$$
Multiplying by $m_{3}$ in numerator and denominator we obtain
$$a_{1}-1=\dfrac{d-w_{0}-w_{1}}{w_{1}}=\dfrac{w_{2}+w_{3}+w_{4}-1}{w_{1}}.$$
Since $w_{2}+w_{3}+w_{4}-1=m_{2}(v_{2}+v_{3}+v_{4})-1,$ we have $\gcd(m_{2},w_{2}+w_{3}+w_{4}-1)=1$. Thus  

\begin{align}
\gcd(m_{2},a_{1}-1)=\gcd\left( m_{2},\dfrac{w_{2}+w_{3}+w_{4}-1}{w_{1}}\right) =1.
\end{align}
Hence, we obtain 
\begin{equation}
\gcd(\Tilde{d},\Tilde{w}_{0})=m_{3}v_{1}.
\end{equation}
On the other hand, since $m_{2}-1=a_{1}v_{1}$, we have
\begin{align}
\gcd(\Tilde{d},\Tilde{w}_{1})&=\gcd(m_{3}m_{2}(m_{2}-1),m_{3}m_{2}v_{1})\nonumber \\
&=m_{3}m_{2}\gcd(m_{2}-1,v_{1})\\
&=m_{3}m_{2}v_{1}\nonumber.
\end{align}
Now we compute $\gcd(\Tilde{d},\Tilde{w}_{2}).$ Here one has 
$$\gcd(\Tilde{d},\Tilde{w}_{2})=m_{2}(m_{2}-1)\gcd(m_{3},1-a_{4}+a_{3}a_{4}).$$
Let us see that $\gcd(m_{3},1-a_{4}+a_{3}a_{4})=1$. Consider $k=\gcd(m_{3},1-a_{4}+a_{3}a_{4})$. By (3.4.7) we have that 
$m_{3}=1+a_{2}a_{3}a_{4}$, then $$\gcd(m_{3},a_{2})=\gcd(m_{3},a_{3})=\gcd(m_{3},a_{4})=1.$$
From the fact that  $m_{3}$ and $1-a_{4}+a_{3}a_{4}$ are divisible by $k$, it follows that  $k$ divides $(1+a_{2}a_{3}a_{4}-(1-a_{4}+a_{3}a_{4}))$. Since  
$\gcd(m_{3},a_{4})=1$, $k$ divides $(1-a_{3}+a_{2}a_{3})$. Repeating the steps given above, one concludes that  $k$ divides 
$(1-a_{2}+a_{2}a_{4})$ as well. On the other hand, from Equation (3.4.4) one has  
$$w_{2}+w_{3}+w_{4}=m_{2}(1-a_{4}+a_{2}a_{4}+1-a_{3}+a_{3}a_{4}+1-a_{2}+a_{2}a_{3}).$$
So $w_{2}+w_{3}+w_{4}$ is divided by $k$. Moreover, we have that $w_{0}=m_{3}$, $w_{1}=m_{3}v_{1}$ and $d+1=w_{0}+w_{1}+w_{2}+w_{3}+w_{4}$, then $k$ divides $(d+1)$. However, since $k$ divides $d$, thus $k=1$. Analogously, it follows 
$\gcd(m_{3},1-a_{3}+a_{2}a_{3})=1$ and $\gcd(1-a_{2}+a_{4}a_{2})=1$. Then, we have     
\begin{equation}   \gcd(\Tilde{d},\Tilde{w}_{2})=\gcd(\Tilde{d},\Tilde{w}_{3})=\gcd(\Tilde{d},\Tilde{w}_{4})=m_{2}(m_{2}-1).
\end{equation}
From (3.4.11), (3.4.12) and (3.4.13), one obtains the numbers $\Tilde{u}_{i}$ and $\Tilde{v}_{i}$ described in (3.4.9):
\begin{equation}
    \Tilde{u}_{0}=m_{2}a_{1}, \ \Tilde{u}_{1}=a_{1}, \ \Tilde{u}_{2}=\Tilde{u}_{3}=\Tilde{u}_{4}=m_{3}.
\end{equation}
and
\begin{equation}
    \Tilde{v}_{0}=a_{1}-1, \  \Tilde{v}_{1}=1, \ \Tilde{v}_{2}=1-a_{4}+a_{3}a_{4}, \ \Tilde{v}_{3}=1-a_{2}+a_{4}a_{2}, \ \Tilde{v}_{4}=1-a_{3}+a_{2}a_{3}.
\end{equation}
Now, we compute $b_{3}$: write  $b_{3}=-S_{0}+S_{1}-S_{2}+S_{3}-S_{4}+S_{5}$, where  
{\small{
\begin{equation*}
    S_{0}=1, \
    S_{1}=\sum \dfrac{1}{\Tilde{v}_{i}}, \ 
    S_{2}=\sum \dfrac{\Tilde{u}_{i_{1}}\Tilde{u}_{i_{2}}}{\Tilde{v}_{i_{1}}\Tilde{v}_{i_{2}}\lcm(\Tilde{u}_{i_{1}},\Tilde{u}_{i_{2}})}, \ 
    S_{3}=\sum \dfrac{\Tilde{u}_{i_{1}}\Tilde{u}_{i_{2}}\Tilde{u}_{i_{3}}}{\Tilde{v}_{i_{1}}\Tilde{v}_{i_{2}}\Tilde{v}_{i_{3}}\lcm(\Tilde{u}_{i_{1}},\Tilde{u}_{i_{2}},\Tilde{u}_{i_{3}})}, 
\end{equation*}
\begin{equation*}
\ S_{4}=\sum \dfrac{\Tilde{u}_{i_{1}}\Tilde{u}_{i_{2}}\Tilde{u}_{i_{3}}\Tilde{u}_{i_{4}}}{\Tilde{v}_{i_{1}}\Tilde{v}_{i_{2}}\Tilde{v}_{i_{3}}\Tilde{v}_{i_{4}}\lcm(\Tilde{u}_{i_{1}},\Tilde{u}_{i_{2}},\Tilde{u}_{i_{3}},\Tilde{u}_{i_{4}})} \ \ \mbox{ and } \ \ \ 
S_{5}=\sum \dfrac{\Tilde{u}_{0}\Tilde{u}_{1}\Tilde{u}_{2}\Tilde{u}_{3}\Tilde{u}_{4}}{\Tilde{v}_{0}\Tilde{v}_{1}\Tilde{v}_{2}\Tilde{v}_{3}\Tilde{v}_{4}\lcm(\Tilde{u}_{0},\Tilde{u}_{1},\Tilde{u}_{2},\Tilde{u}_{3},\Tilde{u}_{4})}.
\end{equation*}
}}

Replacing (3.4.14) and (3.4.15) in the equations given above, after simplifications we obtain:
\begin{align*}
    S_{1} & = \dfrac{a_{1}}{a_{1}-1}+\dfrac{1}{\Tilde{v}_{2}}+\dfrac{1}{\Tilde{v}_{3}}+\dfrac{1}{\Tilde{v}_{4}}. \\
    S_{2} & = \dfrac{a_{1}}{a_{1}-1}+\dfrac{a_{1}^{2}m_{3}}{(a_{1}-1)\lcm(a_{1},m_{3})}\left( \dfrac{1}{\Tilde{v}_{2}}+\dfrac{1}{\Tilde{v}_{3}}+\dfrac{1}{\Tilde{v}_{4}}\right)+m_{3}\left( \dfrac{1}{\Tilde{v}_{2}\Tilde{v}_{3}}+\dfrac{1}{\Tilde{v}_{3}\Tilde{v}_{4}}+\dfrac{1}{\Tilde{v}_{4}\Tilde{v}_{2}}\right). \\
    S_{3} & = \dfrac{a_{1}^{2}m_{3}}{(a_{1}-1)\lcm(a_{1},m_{3})}\left( \dfrac{1}{\Tilde{v}_{2}}+\dfrac{1}{\Tilde{v}_{3}}+\dfrac{1}{\Tilde{v}_{4}}\right) \\
    &\quad + 
    \dfrac{a_{1}^{2}m_{3}^{2}}{(a_{1}-1)\lcm(a_{1},m_{3})}\left( \dfrac{1}{\Tilde{v}_{2}\Tilde{v}_{3}}+\dfrac{1}{\Tilde{v}_{3}\Tilde{v}_{4}}
    +\dfrac{1}{\Tilde{v}_{4}\Tilde{v}_{2}}\right) + \dfrac{m_{3}^{2}}{\Tilde{v}_{2}\Tilde{v}_{3}\Tilde{v}_{4}}. \\
S_{4} & = \dfrac{a_{1}^{2}m_{3}^{2}}{(a_{1}-1)\lcm(a_{1},m_{3})}\left( \dfrac{1}{\Tilde{v}_{2}\Tilde{v}_{3}}+\dfrac{1}{\Tilde{v}_{3}\Tilde{v}_{4}}+\dfrac{1}{\Tilde{v}_{4}\Tilde{v}_{2}}\right) + \dfrac{a_{1}^{2}m_{3}^{3}}{\lcm(a_{1},m_{3})(a_{1}-1)\Tilde{v}_{2}\Tilde{v}_{3}\Tilde{v}_{4}}.\\
    S_{5} & = \dfrac{a_{1}^{2}m_{3}^{3}}{\lcm(a_{1},m_{3})(a_{1}-1)\Tilde{v}_{2}\Tilde{v}_{3}\Tilde{v}_{4}}.
\end{align*}
Replacing these equalities in $b_{3}$, we have
\begin{align*}
    b_{3} & = -1+\left(\dfrac{1}{\Tilde{v_{2}}}+\dfrac{1}{\Tilde{v_{3}}}+\dfrac{1}{\Tilde{v_{4}}}\right)-m_{3}\left(\dfrac{1}{\Tilde{v_{2}}\Tilde{v_{3}}}+\dfrac{1}{\Tilde{v_{3}}\Tilde{v_{4}}}+\dfrac{1}{\Tilde{v_{4}}\Tilde{v_{2}}}\right)+\dfrac{\Tilde{m_{3}}^{2}}{\Tilde{v_{2}}\Tilde{v_{3}}\Tilde{v_{4}}} \\
    & = -1+\dfrac{1}{m_{3}}\left[ \dfrac{(m_{3}-\Tilde{v}_{2})(m_{3}-\Tilde{v}_{3})(m_{3}-\Tilde{v}_{4})}{\Tilde{v}_{2}\Tilde{v}_{3}\Tilde{v}_{4}} +1 \right]. 
\end{align*}

From the matrix equation (3.4.1), we have 
\begin{align*}
    a_{2} & = \dfrac{\Tilde{d}-\Tilde{w}_{3}}{\Tilde{w}_{2}}  = \dfrac{m_{3}-\Tilde{v}_{3}}{\Tilde{v}_{2}},\\
    a_{3} & = \dfrac{\Tilde{d}-\Tilde{w}_{4}}{\Tilde{w}_{3}}  = \dfrac{m_{3}-\Tilde{v}_{4}}{\Tilde{v}_{3}},\\
    a_{4} & = \dfrac{\Tilde{d}-\Tilde{w}_{2}}{\Tilde{w}_{4}}  = \dfrac{m_{3}-\Tilde{v}_{2}}{\Tilde{v}_{4}}.
\end{align*}
Replacing in $b_{3}$, we obtain
$$b_{3}=-1+\dfrac{a_{2}a_{3}a_{4}+1}{m_{3}}.$$
Finally, since $m_{3}=a_{2}a_{3}a_{4}+1$ (Eq. 3.4.9), we conclude that $b_{3}=0.$ Notice that in particular 
\begin{equation}
\left(\dfrac{1}{\Tilde{v_{2}}}+\dfrac{1}{\Tilde{v_{3}}}+\dfrac{1}{\Tilde{v_{4}}}\right)-m_{3}\left(\dfrac{1}{\Tilde{v_{2}}\Tilde{v_{3}}}+\dfrac{1}{\Tilde{v_{3}}\Tilde{v_{4}}}+\dfrac{1}{\Tilde{v_{4}}\Tilde{v_{2}}}\right)+\dfrac{\Tilde{m_{3}}^{2}}{\Tilde{v_{2}}\Tilde{v_{3}}\Tilde{v_{4}}} =1.
\end{equation}
\medskip

\noindent\textbf{Third step:} Now we compute the torsion for this case.  Recall that  by Orlik's algorithm, given as Conjecture 1.1,  we have 
    $H_{3}(L_{f},\mathbb{Z})_{tor} = \mathbb{Z}_{d_{1}}\oplus\cdots\oplus\mathbb{Z}_{d_{r}}$
    where $d_{j}=\prod_{k_{i_{1},\dots, i_{s}}\geq j}c_{i_{1},\dots, i_{s}}.$ Since $k_{i_{1}}=k_{i_{1}i_{2}i_{3}}=0$,  it is sufficient to
    calculate $k_{i_{1}i_{2}}$ and $k_{i_{1}i_{2}i_{3}i_{4}}$, where the associated values $c_{i_{1}i_{2}}$ and $c_{i_{1}i_{2}i_{3}i_{4}}$ are greater than 1.
We obtain $c_{\emptyset}=\gcd(a_{1},m_{3})$, $c_{01}=c_{234}=\dfrac{m_{3}}{\gcd(a_{1},m_{3})}$, $c_{1234}=m_{2},$ and $c_{i_{1},\dots i_{s}}=1$ in other cases. Then, since $k_{\emptyset}=1$ and $k_{234}=0$,  we only need to compute $k_{01}$ and $k_{1234}$. From the definition, we have
\begin{align*}
    k_{01} & = 1-\frac{1}{\Tilde{v}_{0}}-\frac{1}{\Tilde{v_{1}}}+\dfrac{\Tilde{u}_{0}\Tilde{u}_{1}}{\Tilde{v}_{0}\Tilde{v}_{1}\lcm(\Tilde{u}_{0},
\Tilde{u}_{1})}. \\
    \end{align*}

Applying equations (3.4.14) and (3.4.15) we obtain  
   \begin{align*}
k_{01} & = 1-\dfrac{1}{a_{1}-1}-1+\dfrac{m_{2}a_{1}^{2}}{(a_{1}-1)m_{2}a_{1}} \\
    & = 1.
\end{align*}
On the other hand, let us compute $k_{1234}$. We express this number as 
$$k_{1234}=\Sigma_{0}-\Sigma_{1}+\Sigma_{2}-\Sigma_{3}+\Sigma_{4},$$ where

\begin{align*}
    \Sigma_{0} & = 1. \\
    \Sigma_{1} & = \dfrac{1}{\Tilde{v}_{1}}+\dfrac{1}{\Tilde{v}_{2}}+\dfrac{1}{\Tilde{v}_{3}}+\dfrac{1}{\Tilde{v}_{4}},\\
    \Sigma_{2} & = \dfrac{\Tilde{u}_{1}\Tilde{u}_{2}}{\Tilde{v}_{1}\Tilde{v}_{2}\lcm(\Tilde{u}_{1},\Tilde{u}_{2})} +
   \dfrac{\Tilde{u}_{1}\Tilde{u}_{3}}{\Tilde{v}_{1}\Tilde{v}_{3}\lcm(\Tilde{u}_{1},\Tilde{u}_{3})} +  \dfrac{\Tilde{u}_{1}\Tilde{u}_{4}}{\Tilde{v}_{1}\Tilde{v}_{4}\lcm(\Tilde{u}_{1},\Tilde{u}_{4})} \\
    & \quad +  \dfrac{\Tilde{u}_{2}\Tilde{u}_{3}}{\Tilde{v}_{2}\Tilde{v}_{3}\lcm(\Tilde{u}_{2},\Tilde{u}_{3})} +  \dfrac{\Tilde{u}_{2}\Tilde{u}_{4}}{\Tilde{v}_{2}\Tilde{v}_{4}\lcm(\Tilde{u}_{2},\Tilde{u}_{4})} +  \dfrac{\Tilde{u}_{3}\Tilde{u}_{4}}{\Tilde{v}_{3}\Tilde{v}_{4}\lcm(\Tilde{u}_{3},\Tilde{u}_{4})},\\
    \Sigma_{3} & = \dfrac{\Tilde{u}_{1}\Tilde{u}_{2}\Tilde{u}_{3}}{\Tilde{v}_{1}\Tilde{v}_{2}\Tilde{v}_{3}\lcm(\Tilde{u}_{1},\Tilde{u}_{2},\Tilde{u}_{3})}
   + \dfrac{\Tilde{u}_{1}\Tilde{u}_{2}\Tilde{u}_{4}}{\Tilde{v}_{1}\Tilde{v}_{2}\Tilde{v}_{4}\lcm(\Tilde{u}_{1},\Tilde{u}_{2},\Tilde{u}_{4})} \\
  & \quad + \dfrac{\Tilde{u}_{1}\Tilde{u}_{3}\Tilde{u}_{4}}{\Tilde{v}_{1}\Tilde{v}_{3}\Tilde{v}_{4}\lcm(\Tilde{u}_{1},\Tilde{u}_{3},\Tilde{u}_{4})} + \dfrac{\Tilde{u}_{2}\Tilde{u}_{3}\Tilde{u}_{4}}{\Tilde{v}_{2}\Tilde{v}_{3}\Tilde{v}_{4}\lcm(\Tilde{u}_{2},\Tilde{u}_{3},\Tilde{u}_{4})}, \\
    \Sigma_{4} & = \dfrac{\Tilde{u}_{1}\Tilde{u}_{2}\Tilde{u}_{3}\Tilde{u}_{4}}{\Tilde{v}_{1}\Tilde{v}_{2}\Tilde{v}_{3}\Tilde{v}_{4}\lcm(\Tilde{u}_{1},\Tilde{u}_{2},\Tilde{u}_{3},\Tilde{u}_{4})}.
\end{align*}

Bearing in mind equations (2.4.13) and (2.4.14) one obtains
\begin{align*}
    \Sigma_{1} & = 1+\dfrac{1}{\Tilde{v}_{2}}+\dfrac{1}{\Tilde{v}_{3}}+\dfrac{1}{\Tilde{v}_{4}}, \\
    \Sigma_{2} & = \dfrac{a_{1}m_{3}}{\lcm(a_{1},m_{3})}\left( \dfrac{1}{\Tilde{v}_{2}}+\dfrac{1}{\Tilde{v}_{3}}+\dfrac{1}{\Tilde{v}_{4}}\right) + m_{3}\left( \dfrac{1}{\Tilde{v}_{2}\Tilde{v}_{3}}+\dfrac{1}{\Tilde{v}_{2}\Tilde{v}_{4}}+\dfrac{1}{\Tilde{v}_{3}\Tilde{v}_{4}}\right), \\
    \Sigma_{3} & = \dfrac{a_{1}m_{3}^{2}}{\lcm(a_{1},m_{3})}\left( \dfrac{1}{\Tilde{v}_{2}\Tilde{v}_{3}}+\dfrac{1}{\Tilde{v}_{2}\Tilde{v}_{4}}+\dfrac{1}{\Tilde{v}_{3}\Tilde{v}_{4}}\right) + \dfrac{m_{3}^{2}}{\Tilde{v}_{2}\Tilde{v}_{3}\Tilde{v}_{4}},\\
    \Sigma_{4} & = \dfrac{a_{1}m_{3}^{3}}{\lcm(a_{1},m_{3})\Tilde{v}_{2}\Tilde{v}_{3}\Tilde{v}_{4}}.
\end{align*}
Thus,  
$$k_{1234}=\left(\dfrac{a_{1}m_{3}}{\lcm(a_{1},m_{3})}-1\right)\left[ \dfrac{1}{\Tilde{v}_{2}}+\dfrac{1}{\Tilde{v}_{3}}+\dfrac{1}{\Tilde{v}_{4}} - m_{3}\left( \dfrac{1}{\Tilde{v}_{2}\Tilde{v}_{3}}+\dfrac{1}{\Tilde{v}_{2}\Tilde{v}_{4}}+\dfrac{1}{\Tilde{v}_{3}\Tilde{v}_{4}}\right) + \dfrac{m_{3}^{2}}{\Tilde{v}_{2}\Tilde{v}_{3}\Tilde{v}_{4}}\right].$$
By Equation (3.4.15), we have 
$$\dfrac{1}{\Tilde{v}_{2}}+\dfrac{1}{\Tilde{v}_{3}}+\dfrac{1}{\Tilde{v}_{4}} - m_{3}\left( \dfrac{1}{\Tilde{v}_{2}\Tilde{v}_{3}}+\dfrac{1}{\Tilde{v}_{2}\Tilde{v}_{4}}+\dfrac{1}{\Tilde{v}_{3}\Tilde{v}_{4}}\right) + \dfrac{m_{3}^{2}}{\Tilde{v}_{2}\Tilde{v}_{3}\Tilde{v}_{4}}=1.$$
Hence 
\begin{equation}
    k_{1234}=\dfrac{a_{1}m_{3}}{\lcm(a_{1},m_{3})}-1 = \gcd(a_{1},m_{3})-1. 
\end{equation}
Then, there are three possible values for $H_{3}(L_{f^T},\mathbb{Z})$:
\begin{itemize}
    \item If $\gcd(a_{1},m_{3})-1<1$, then $H_{3}(L_{f^T},\mathbb{Z})=\mathbb{Z}_{m_{3}}$.
    \item If $\gcd(a_{1},m_{3})-1=1$, then $H_{3}(L_{f^T},\mathbb{Z})=\mathbb{Z}_{d}$.
    \item If $\gcd(a_{1},m_{3})-1>1$, then $H_{3}(L_{f^T},\mathbb{Z})=\mathbb{Z}_{d}\oplus \underbrace{\mathbb{Z}_{m_{2}}\oplus\cdots\oplus\mathbb{Z}_{m_{2}}}_{(\gcd(a_{1},m_{3})-2)-times}$.
\end{itemize}

\noindent\textbf{Fourth step.} Next we will compute the Milnor number $\tilde{\mu}$: the steps are very similar to the ones taken in the fourth step in Case II, but unlike that case, the first two weights $\tilde{w}_{0}$ and $\tilde{w}_{1}$ are different than the original ones $w_0$ and $w_1.$

 We will use the formula:
\begin{equation}
    \tilde{\mu}=\left(\dfrac{\tilde{d}-\tilde{w}_{0}}{\tilde{w}_{0}}\right)\left(\dfrac{\tilde{d}-\tilde{w}_{1}}{\tilde{w}_{1}}\right)\left(\dfrac{\tilde{d}-\tilde{w}_{2}}{\tilde{w}_{2}}\right)\left(\dfrac{\tilde{d}-\tilde{w}_{3}}{\tilde{w}_{3}}\right)\left(\dfrac{\tilde{d}-\tilde{w}_{4}}{\tilde{w}_{4}}\right).
\end{equation}
Since $a_{2}a_{3}a_{4}=m_{3}-1$ and  
$$a_{2}a_{3}a_{4}=\left(\dfrac{\tilde{d}-\tilde{w}_{3}}{\tilde{w}_{2}}\right)\left(\dfrac{\tilde{d}-\tilde{w}_{4}}{\tilde{w}_{3}}\right)\left(\dfrac{\tilde{d}-\tilde{w}_{2}}{\tilde{w}_{4}}\right),$$
we obtain
\begin{equation}
    \left(\dfrac{\tilde{d}-\tilde{w}_{2}}{\tilde{w}_{2}}\right)\left(\dfrac{\tilde{d}-\tilde{w}_{3}}{\tilde{w}_{3}}\right)\left(\dfrac{\tilde{d}-\tilde{w}_{4}}{\tilde{w}_{4}}\right)=m_{3}-1.
\end{equation}

On the other hand, we have 
\begin{equation}
    \dfrac{\tilde{d}-\tilde{w}_{0}}{\tilde{w}_{0}}=\dfrac{m_{2}a_{1}-a_{1}+1}{a_{1}-1} \ \  \mbox{ and } \dfrac{\tilde{d}-\tilde{w}_{1}}{\tilde{w}_{1}}=a_{1}-1.
\end{equation}
Replacing (3.4.19) and (3.4.20) in (3.4.18), we obtain 
$$\tilde{\mu}=(m_{2}a_{1}-a_{1}+1)(m_{3}-1).$$
As $a_{1}=\frac{m_{2}-1}{v_{1}}$, we finally have 
\begin{equation*}
    \tilde{\mu}=\left( \dfrac{(m_{2}-1)^{2}}{v_{1}}+1\right)(m_{3}-1).
\end{equation*}

With these four cases, we finish the proof of Theorem 3.2.


\begin{remark}
We notice that in  the case of weights $\textbf{w}=(w_{0},w_{1},\dots, w_{4})$ with a well-formed 3-fold $X_{f}$ where $f$ is a chain-cycle polynomials, the hypersurface $X_{f^{T}}$ obtained through the Berglund-Hubsch method does not remain  well-formed. For instance, we can take the vector weight $\textbf{w}=(881,881,465,99,318)$ of the list of Johnson-Koll\'ar. Here we have the hypersurface defined by the chain-cycle polynomial:
$$f=z_{0}^{3}+z_{0}z_{1}^{2}+z_{4}z_{2}^{5}+z_{2}z_{3}^{22}+z_{3}z_{4}^{8}$$
with degree $d=2643$. If we denote by $A$ it matrix of exponents,  through BH-transpose rule, $A^T$ is given by 

\begin{equation*}
    A^{T}=\begin{bmatrix}
      3 & 1 & 0 & 0 & 0 \\
      0 & 2 & 0 & 0 & 0 \\
      0 & 0 & 5 & 1 & 0 \\
      0 & 0 & 0 & 22 & 1 \\
      0 & 0 & 1 & 0 & 8 
    \end{bmatrix}
\end{equation*}
Solving the matrix equation $A^{T}\tilde{W}=D$, where $D$ is the column matrix $5\times 1$ and  $d=2643$, we obtain the vector weights $\mathbf{w}_{T}=(881,2643,1014,216,534)$, which  produces the 3-fold $X_{f_{T}}$ that is not well-formed since  
$\gcd(2643,1014,216,534)\neq 1.$ 
\end{remark}

Actually, we  have the following result:

\begin{theorem} Let $f$ be an invertible polynomial as in Theorem 3.2
\begin{enumerate}
    \item If $f$ lies in either case I, II or II,  then the resulting  hypersurface via BH-transpose rule $X_{f^T}$ is well-formed.
    \item If $f$ lies in case IV, then the resulting  hypersurface via BH-transpose rule $X_{f^T}$ is not well-formed.
\end{enumerate}
\end{theorem}
\noindent{\it Proof.} We first prove the first statement.  

Recall that if $X_d \subset \mathbb{P}(\mathbf{w})$ is a weighted homogeneous hypersurface, then $X_d$ is well-formed if and only if $\mathbb{P}(\mathbf{w})$ is well-formed and for all $0 \leq i<j \leq n$ we have $\operatorname{gcd}\left(w_0, \ldots, \hat{w}_i, \ldots, \hat{w}_j, \ldots, w_n\right) \mid d.$  Thus the weights $\tilde{w}_{i}$'s must  verify the following conditions:
\begin{itemize}
    \item[(i)] $\gcd(\tilde{w}_{0},\dots,\hat{\tilde{w}}_{i},\dots,\tilde{w}_{4})=1,$ 
    \item[(ii)] $\gcd(\tilde{w}_{0},\dots,\hat{\tilde{w}}_{i},\dots,\hat{\tilde{w}}_{j},\dots,\tilde{w}_{4}) | d.$
\end{itemize}
For  BP-cycle singularities consider  the matrix equation
$A^{T}\Tilde{W}=D.$ We obtain 
\begin{align*}
       a_{0}\tilde{w}_{0} & = d,\\
       a_{1}\tilde{w}_{1} & = d,\\ 
       a_{2}\tilde{w}_{2}+\tilde{w}_{3}& =d,\\ 
       a_{3}\tilde{w}_{3}+\tilde{w}_{4} & =d,\\ 
       a_{4}\tilde{w}_{4}+\tilde{w}_{2}& =d.   
   \end{align*}
Let us show  condition (i): it is  sufficient to study the following two cases:
    \begin{itemize}
        \item First, if $\gcd(\tilde{w}_{0},\tilde{w}_{1},\tilde{w}_{2},\tilde{w}_{3})=r$, then we have $r | (d+1-\tilde{w}_{4})$. As $\tilde{w}_{4}=d-a_{3}\tilde{w}_{3}$, we obtain $r | (1+a_{3}\tilde{w}_{3})$ and $r | \tilde{w}_{3}$. Thus $r=1$.
        
        \item On the other hand, if $\gcd(\tilde{w}_{0},\tilde{w}_{2},\tilde{w}_{3},\tilde{w}_{4})=r$, then we have $r | \tilde{w}_{1}$ (since 
        $\tilde{w}_{0}=\tilde{w}_{1}$ in  BP-cycle polynomials). Then $r | (d+1)$, which implies $r=1$.
    \end{itemize}
    
To show  condition (ii), we will assume, without loss of generality,  two cases:
\begin{itemize}
    \item If $r=\gcd(\tilde{w}_{2},\tilde{w}_{3},\tilde{w}_{4})$, then  $r | a_{2}\tilde{w}_{2}+\tilde{w}_{3}$. Therefore $r | d$.
    \item  $r=\gcd(\tilde{w}_{i},\tilde{w}_{j},\tilde{w}_{j})$, where some $i,j$ or $k$ is 0 or 1. In this case, for instance, if $i=0$, then $r | \tilde{w}_{0}$ and $a_{0}\tilde{w}_{0}=d$. Therefore $r | d$. 
\end{itemize}
For  cycle-cycle singularities, we obtain weights $\tilde{w}_{i}$'s solving the matrix equation $A^{T}\Tilde{W}=D$, such that 
\begin{align*}
       a_{0}\tilde{w}_{0}+\tilde{w}_{1}& = d,\\
       a_{1}\tilde{w}_{1}+\tilde{w}_{0}& =d,\\ 
       a_{2}\tilde{w}_{2}+\tilde{w}_{3}& =d,\\ 
       a_{3}\tilde{w}_{3}+\tilde{w}_{4}& =d,\\ 
       a_{4}\tilde{w}_{4}+\tilde{w}_{2}& =d.   
   \end{align*}
We study two cases, in a similar way as above. To show  condition (i), consider two cases: 
\begin{itemize}
        \item First, if $\gcd(\tilde{w}_{0},\tilde{w}_{1},\tilde{w}_{2},\tilde{w}_{3})=r$, then we have $r | (d+1-\tilde{w}_{4})$. As $\tilde{w}_{4}=d-a_{3}\tilde{w}_{3}$, we obtain $r | (1+a_{3}\tilde{w}_{3})$ and $r | \tilde{w}_{3}$. Thus $r=1$.
        
        \item On the other hand, if $\gcd(\tilde{w}_{0},\tilde{w}_{2},\tilde{w}_{3},\tilde{w}_{4})=r$, then we have $r | (d+1-\tilde{w}_{1})$ . As $\tilde{w}_{1}=d-a_{0}\tilde{w}_{0}$, we obtain $r | (1+a_{0}\tilde{w}_{0})$. since $r | \tilde{w}_{0}$, we conclude that $r=1$.
    \end{itemize}
Condition (ii) is proven as follows: let  $r=\gcd(\tilde{w}_{i},\tilde{w}_{j},\tilde{w}_{k}) | d$, we can consider the two following situations:
\begin{itemize}
    \item If $\{ 0,1\}\subset \{ i,j,k\}$. In these case, as $a_{0}\tilde{w}_{0}+\tilde{w}_{1}= d$, then $r | d$.
    \item If $\{ 0,1\} \not\subset\{i,j,k\}$, we can suppose without loss of generality that $i=2$ and $j=3$. Then $r | d$.
\end{itemize}
 For singularities of  cycle type, we have that the weights $\tilde{w}_{i}$'s that solve  $A^{T}\Tilde{W}=D$ satisfy the equations \begin{equation*}a_{0}\tilde{w}_{0}+\tilde{w}_{1}=d, \ \ a_{0}\tilde{w}_{1}+\tilde{w}_{2}=d, \ \ a_{2}\tilde{w}_{2}+\tilde{w}_{3}=d, \ \ a_{3}\tilde{w}_{3}+\tilde{w}_{4}=d, \ \ a_{4}\tilde{w}_{4}+\tilde{w}_{0}=d.
\end{equation*}
In order to verify  condition (i), we can consider without loss of generality that $$\gcd(\tilde{w}_{0},\tilde{w}_{1},\tilde{w}_{2},\tilde{w}_{3})=r,$$ which implies that $r | (d+1-\tilde{w}_{4})$. As $r | \tilde{w}_{3}$, we have that $r | (d-\tilde{w}_{4})$. Thus $r=1$. 
To show condition (ii), let  $r=\gcd(\tilde{w}_{i},\tilde{w}_{j},\tilde{w}_{k}) | d$, we can consider two situations:
 If $i$ and $j$ are consecutive numbers, we have $j=i+1$. Then $r | a_{i}\tilde{w}_{i}$ and $r | \tilde{w}_{i+1}$, thus $r | d$.
 In another case: $i=0$, $j=2$ and $k=4$. Here $r | a_{4}\tilde{w}_{4}$ and $r | \tilde{w}_{0}$, hence $r | d$.

Now, we  prove item $(2)$: recall that the vector weight $\textbf{w}_{T}=(\Tilde{w}_{0},\Tilde{w}_{1},\Tilde{w}_{2},\Tilde{w}_{3},\Tilde{w}_{4})$ is calculated  solving the following matrix equation:
\begin{equation*}
A^{T}\Tilde{W}=\Tilde{D}
\end{equation*}
where $\Tilde{W}=(\Tilde{w}_{0},\Tilde{w}_{1},\Tilde{w}_{2},\Tilde{w}_{3},\Tilde{w}_{4})^{T}$ and $\Tilde{D}=(\Tilde{d},\Tilde{d},\Tilde{d},\Tilde{d},\Tilde{d})^{T}$.

\vspace{0.2cm}

This matrix equation is equivalent to
 \begin{equation*}
        \begin{bmatrix}
          a_{0} & 1 & 0 & 0 & 0 & -1 \\
      0 & a_{1} & 0 & 0 & 0 & -1\\
      0 & 0 & a_{2} & 1 & 0 & -1\\
      0 & 0 & 0 & a_{3} & 1 & -1\\
      0 &  0 & 1 & 0 & a_{4} & -1  
        \end{bmatrix} \begin{bmatrix}
            \Tilde{w}_{0} \\
            \Tilde{w}_{1} \\
            \Tilde{w}_{2} \\
            \Tilde{w}_{3} \\
            \Tilde{w}_{4} \\
            \Tilde{d}
        \end{bmatrix} = \begin{bmatrix}
            0 \\ 
            0 \\
            0 \\
            0 \\
            0
        \end{bmatrix}
    \end{equation*}
As we have seen, the dimension of the solution space of this equations system is 1 and   a nonzero vector solution with integer entries is 

$$\Tilde{w_{1}}=\dfrac{\Tilde{d}}{a_{1}}=m_{3}m_{2}v_{1}$$
$$\Tilde{w_{0}}=\dfrac{\Tilde{d}-\Tilde{w_{1}}}{a_{0}}=\dfrac{m_{3}m_{2}(m_{2}-1)-m_{3}m_{2}v_{1}}{m_{2}}=m_{3}v_{1}(a_{1}-1)$$
$$\Tilde{w}_{2}=\dfrac{\Tilde{d}(1-a_{4}+a_{4}a_{3})}{1+a_{2}a_{3}a_{4}}=m_{2}(m_{2}-1)(1-a_{4}+a_{3}a_{4})$$
$$\Tilde{w}_{3}=\dfrac{\Tilde{d}(1-a_{2}+a_{2}a_{4})}{1+a_{2}a_{3}a_{4}}=m_{2}(m_{2}-1)(1-a_{2}+a_{2}a_{4})$$
$$\Tilde{w}_{4}=\dfrac{\Tilde{d}(1-a_{3}+a_{2}a_{3})}{1+a_{2}a_{3}a_{4}}=m_{2}(m_{2}-1)(1-a_{3}+a_{2}a_{3})$$

Notice that $\tilde{w}_{1},\tilde{w}_{2},\tilde{w}_{3}$ and $\tilde{w}_{4}$ have common factor $m_{2}$. Now we will show that $\gcd(\tilde{w}_{1},m_{2})=1$. As $\gcd(m_{3},m_{2})=1$ and $m_{2}=a_{1}v_{1}+1$, then $\gcd(\tilde{w}_{1},m_{2})=\gcd(v_{1}(a_{1}-1),m_{2})=\gcd(a_{1}-1,m_{2})$. From equation (3.4.10) we have that $\gcd(a_{1}-1,m_{2})=1$, therefore $\gcd(\tilde{w}_{1},m_{2})=1$. That is,  any solution vector $(\tilde{w}_{0},\tilde{w}_{1},\tilde{w}_{2},\tilde{w}_{3},\tilde{w}_{4})$ with integer entries always contains $m_{2}$ in all but one weight. Thus,  
$X_{f^T}$ is not well-formed.

\hfill$\square$

\begin{remark} In \cite{BGN2}, the authors emphasized the existence of  septuplets in their table of rational homology Sasaki-Einstein 7-spheres: a set of seven twins, not related at first glance. From the perspective of  Theorem 3.1 (and  due to the fact that each of the elements there admits a unique representation as a cycle polynomial) it follows that this set actually consists of octuplets since each of these are related by pairs, that is, there is a set of four elements generating the octuplets: each generator producing its BH transpose dual. However one of these dual partners is missing in the table given in \cite{BGN2}, the one corresponding to the data  $(157, 545, 1051, 1401, 2608)$  and degree $d=5761$, given through BH rule by $(148, 477,1871,1321,1945)$ with degree $d=5761.$ This element is exhibited in Table 1 in \cite{CL} that listed 52 rational homology spheres admitting Sasaki-Einstein metrics.
\end{remark}

We have the following result for links lying in cases I, II or III, the proof is an immediate consequence of Theorem 3.1 

\begin{corollary} If $X_f$ is an hypersurface whose corresponding link lies in Case I, II or II, then the BH-transpose rule determines a hypersurface $X_{f^T}$ whose associated link $L_{f^T}$ is a twin for the link $L_f,$ provided $f$ is of cycle type, cycle-cycle type or BP-cycle type.  Pictorially, we have the  diagram

\begin{center}
\begin{tikzcd}
L_{f} \arrow{r}{twin} \arrow{d} & L_{f^T}\arrow{d}\\
{X}_{f} \arrow{r}{BH} & X_{f^T}.\\
\end{tikzcd}
\end{center}
\end{corollary}

\hfill$\square$

For links lying in case IV, we have the next corollary:
\begin{corollary}
    Let us consider the links  $L_f$ and $L_g$ as in case IV, that are twins, 
      then the links $L_{f^T}$ and  $L_{g^T}$,  produced by the BH-transpose rule are twins. Pictorially, we have the diagram
\begin{center}
\begin{tikzcd}
L_{f} \arrow{r}{BH} \arrow[leftrightarrow]{d}[swap]{twin} & L_{f^T}\\
L_{g} \arrow{r}{BH} &L_{g^T}\arrow[leftrightarrow]{u}[swap]{twin}.\\
\end{tikzcd}
\end{center}

\end{corollary}

\proof We consider a pair of twins $\textbf{w}=(w_{0},w_{1},w_{2},w_{3},w_{4})$ and $\textbf{w}'=(w_{0}',w_{1}',w_{2}',w_{3}',w_{4}')$ with hypersurfaces $X_{f}$ and $X_{g}$, where $f$ and $g$ are chain-cycle polynomials and links associated $L_{f}$ and $L_{g}$ which has the same torsion $H_{3}=(\mathbb{Z}_{m_{3}})^{\alpha(w)+1}$. By the BH-transpose rule, we obtain, respectively, the weight vectors $\tilde{\textbf{w}}$ and   $\tilde{\textbf{w}}'$ which are not well-formed and whose torsion and Milnor number are given by formulas:
\begin{itemize}
    \item If $\gcd(a_{1},m_{3})-1<1$, then $H_{3}(L_{f},\mathbb{Z})=\mathbb{Z}_{m_{3}}$.
    \item If $\gcd(a_{1},m_{3})-1=1$, then $H_{3}(L_{f},\mathbb{Z})=\mathbb{Z}_{d}$.
    \item If $\gcd(a_{1},m_{3})-1>1$, then $H_{3}(L_{f},\mathbb{Z})=\mathbb{Z}_{d}\oplus \underbrace{\mathbb{Z}_{m_{2}}\oplus\cdots\oplus\mathbb{Z}_{m_{2}}}_{(\gcd(a_{1},m_{3})-2)-times}$.
\end{itemize}

and 

\begin{equation*}
    \tilde{\mu}=\left( \dfrac{(m_{2}-1)^{2}}{v_{1}}+1\right)(m_{3}-1).
\end{equation*}

Since twins have same $m_{3}$ and $m_{2}$, we obtain that the torsion and the Milnor number of the links of  $\tilde{\textbf{w}}$ and $\tilde{\textbf{w}}'$ are equals.
So, it remains to show that $\tilde{\textbf{w}}$ and $\tilde{\textbf{w}}'$ have the same degree. Since we are working with polynomials as in case IV, we can write the weight vectors $\textbf{w}$ and $\textbf{w}'$  as:
$$\textbf{w}=(w_{0},w_{1},w_{2},w_{3},w_{4})=(m_{3},m_{3}v_{1},m_{2}v_{2},m_{2}v_{3},m_{2}v_{4})$$
and 
$$\textbf{w}'=(w_{0}',w_{1}',w_{2}',w_{3}',w_{4}')=(m_{3},m_{3}v_{1},m_{2}v_{2}',m_{2}v_{3}',m_{2}v_{4}').$$

From the BH-transpose rule, we obtain the associated vector weights $\tilde{\textbf{w}}$ and $\tilde{\textbf{w}}'$:
$$\tilde{\textbf{w}}=(m_{3}m_{2}v_{1},m_{3}v_{1}(a_{1}-1),m_{2}(m_{2}-1)v_{2}',m_{2}(m_{2}-1)v_{3}',m_{2}(m_{2}-1)v_{4}')$$
and 
$$\tilde{\textbf{w}}'=(m_{3}m_{2}v_{1},m_{3}v_{1}(a_{1}-1),m_{2}(m_{2}-1)v_{2},m_{2}(m_{2}-1)v_{3},m_{2}(m_{2}-1)v_{4})$$
both with degree $\tilde{d}=m_{3}m_{2}(m_{2}-1)$.
Moreover, as $m_{2}=a_{1}v_{1}+1$, we can simplify the expressions given above: 
$$\tilde{\textbf{w}}=(m_{3}m_{2}, m_{3}(a_{1}-1), m_{2}a_{1}v_{2}', m_{2}a_{1}v_{3}',m_{2}a_{1}v_{4}')$$
and
$$\tilde{\textbf{w}}'=(m_{3}m_{2},m_{3}(a_{1}-1),m_{2}a_{1}v_{2},m_{2}a_{1}v_{3},m_{2}a_{1}v_{4})$$
with degree $\tilde{d}=m_{3}m_{2}a_{1}$. Since $\gcd(m_{3},m_{2})=1$, $\gcd(m_{2},a_{1}-1)=1$, $\gcd(m_{3},v_{i})=1$ and $\gcd(m_{3},v_{i}')=1$ for $i=1,2,3$, we have that all entries of $\tilde{\textbf{w}}$ and $\tilde{\textbf{w}}'$ have as common factor  $\gcd(m_{3},a_{1})$. Finally, as $$v_{2}+v_{3}+v_{4}=v_{2}'+v_{3}'+v_{4}',$$
we conclude that the degrees are the same.
\hfill$\square$
\medskip


\section{Sasaki-Einstein metrics on links under BH-transpose rule} 
Let us consider the lists of links that are  rational homology 7-spheres admitting Sasaki-Einstein metrics given in \cite{BGN2} and \cite{CL}.  Now, if one applies the BH-transpose rule to members of any of these  two lists of links that come from hypersurface singularities of cycle type, BP-cycle type or cycle-cycle type (which fall in one the three fist cases in the proof of Theorem 2.1) the twin ends up being a  member of one those lists, and hence it admits Sasaki-Einstein metrics.
 Unlike  these three cases, Case IV  manufactures non-twins, thus we will focus on this sort of invertible polynomial.  The theorem presented in this section show that for links of chain-cycle polynomials, the  BH-transpose rule also  preserves Sasaki-Einstein metrics provided certain plausible conditions suggested by the data found in \cite{BGN2} and in \cite{CL} are satisfied. 
\medskip

First we have the following lemma:

\begin{lemma}   Consider the  weight vector $\textbf{w}=(w_{0},w_{1},w_{2},w_{3},w_{4})$ that can be written as
\begin{equation}
\mathbf{w}=(m_{3}v_{0},m_{3}v_{1},m_{2}v_{2},m_{2}v_{3},m_{2}v_{4})
\end{equation}
and a hypersurface $X_{f}$ defined by the polynomial type chain-cycle of degree $d$:
$$f=z_{0}^{a_{0}}+z_{0}z_{1}^{a_{1}}+z_{4}z_{2}^{a_{2}}+z_{2}z_{3}^{a_{3}}+z_{3}z_{4}^{a_{4}}$$
where $d=m_{3}m_{2}$, with $\gcd(m_{3},m_{2})=1$ and index $I=1$.
If the link $L_{f}$  associated to $f$ is a rational homological sphere and the exponents $a_{2}$, $a_{3}$ and $a_{4}$ are different, then $L_{f}$ admits a twin.
\end{lemma}   
\noindent{\it Proof.} Notice that components of $\textbf{w}$ satisfy the matrix equation

\begin{equation}
        \begin{bmatrix}
          a_{0} & 0 & 0 & 0 & 0 \\
      1 & a_{1} & 0 & 0 & 0 \\
      0 & 0 & a_{2} & 0 & 1 \\
      0 & 0 & 1 & a_{3} & 0 \\
      0 &  0 & 0 & 1 & a_{4}   
        \end{bmatrix} \begin{bmatrix}
            w_{0} \\
            w_{1} \\
            w_{2} \\
            w_{3} \\
            w_{4}
        \end{bmatrix} = \begin{bmatrix}
            d \\ 
            d \\
            d \\
            d \\
            d 
        \end{bmatrix}
    \end{equation}

Indeed, we can define the polynomial of  chain-cycle type
$$g=z_{0}^{a_{0}}+z_{0}z_{1}^{a_{1}}+z_{4}z_{2}^{a_{3}}+z_{2}z_{3}^{a_{2}}+z_{3}z_{4}^{a_{4}}.$$
Next, we will look for weights $\textbf{w}'=(w_{0}',w_{1}',w_{2}',w_{3}',w_{4}')$ such that these  solve the matrix equation 
\begin{equation}
        \begin{bmatrix}
          a_{0} & 0 & 0 & 0 & 0 \\
      1 & a_{1} & 0 & 0 & 0 \\
      0 & 0 & a_{3} & 0 & 1 \\
      0 & 0 & 1 & a_{2} & 0 \\
      0 &  0 & 0 & 1 & a_{4}   
        \end{bmatrix} \begin{bmatrix}
            w_{0}' \\
            w_{1}' \\
            w_{2}' \\
            w_{3}' \\
            w_{4}'
        \end{bmatrix} = \begin{bmatrix}
            d \\ 
            d \\
            d \\
            d \\
            d 
        \end{bmatrix}
    \end{equation}
Since the matrix  above is invertible, the equation has unique solution: $w_{0}'=w_{0}, w_{1}=w_{1}'$ and 
$$ w_{2}'=\dfrac{d}{a_{2}a_{3}a_{4}+1}(a_{2}a_{4}-a_{2}+1), \  w_{3}'=\dfrac{d}{a_{2}a_{3}a_{4}+1}(a_{3}a_{4}-a_{4}+1),  \ w_{4}'=\dfrac{d}{a_{2}a_{3}a_{4}+1}(a_{3}a_{2}-a_{3}+1).$$
As $L_{f}$ is a rational homological sphere, we have that $a_{2}a_{3}a_{4} + 1 = m_{3} $.  Thus, these weights are integers. Furthermore, we can verify that
$$w_{0}'+w_{1}'+w_{2}'+w_{3}'+w_{4}'=d+1$$
and $X_{f}$ is well-formed in a similar way as in the case of polynomials of cycle-cycle type. If we put $v_{2}'=a_{2}a_{4}-a_{2}+1$, $v_{3}'=a_{3}a_{4}-a_{4}+1$ and $v_{4}'=a_{3}a_{2}-a_{3}+1$, then we can write the weights as
\begin{equation}
    w_{2}'=m_{2}v_{2}', \ \ w_{3}'=m_{2}v_{3}', \ \ w_{4}'=m_{2}v_{4}'.
\end{equation}
Thus, taking into account that $w_{0}=w_{0}'$ and $w_{1}=w_{1}'$, we can express the vector weights $\textbf{w}'$:
$$\textbf{w}'=(m_{3}v_{0},m_{3}v_{1},m_{2}v_{2}',m_{2}v_{3}',m_{2}v_{4}').$$
Next, we will show that $b_{3}'=0$ and calculate its torsion. Due to the remark given in B), it is sufficient to prove that $\beta(\mathbf{w}')=1$. From (4.0.1) and (4.0.3) we have
\begin{equation}
    a_{2} = \dfrac{d-w_{4}}{w_{2}} = \dfrac{d-w_{2}'}{w_{3}'} \ \ \
    a_{3} = \dfrac{d-w_{2}}{w_{3}} = \dfrac{d-w_{4}'}{w_{2}'} \ \ \ 
    a_{4} = \dfrac{d-w_{3}}{w_{4}} = \dfrac{d-w_{3}'}{w_{4}'}.
\end{equation}
Replacing (4.0.1) and (4.0.4) in the formula above, and multiplying, we obtain:
\begin{equation*}
    a_{2}a_{3}a_{4}=\dfrac{(m_{3}-v_{2})(m_{3}-v_{3})(m_{3}-v_{4})}{v_{2}v_{3}v_{4}}=\dfrac{(m_{3}-v_{2}')(m_{3}-v_{3}')(m_{3}-v_{4}')}{v_{2}v_{3}v_{4}'}.
\end{equation*}
Therefore, we have that 
   $$\beta(\mathbf{w})= \dfrac{m_{3}^{2}-(v_{2}+v_{3}+v_{4})m_{3}+v_{2}v_{3}+v_{2}v_{4}+v_{3}v_{4}}{v_{2}v_{3}v_{4}}$$  is equal to the expression   
   $$\dfrac{m_{3}^{2}-(\tilde{v}_{2}+\tilde{v}_{3}+\tilde{v}_{4})m_{3}+\tilde{v}_{2}\tilde{v}_{3}+\tilde{v}_{2}\tilde{v}_{4}+\tilde{v}_{3}\tilde{v}_{4}}{\tilde{v}_{2}\tilde{v}
    _{3}\tilde{v}_{4}}=\beta(\mathbf{w}').$$
     
As $\beta(\mathbf{w})=1$, then $\beta(\mathbf{w}')=1$. Thus $b_{3}'=0$.

Now we prove that the links $L_{f}$ and $L_{g}$ have the same torsion. From the remark given in B) after formula $2.2.8$, we know   
$$H_{3}(L_{f},\mathbb{Z})=(\mathbb{Z}_{m_{3}})^{\alpha(\textbf{w})+1}$$
where 
$$\alpha(\textbf{w})+1=\dfrac{m_{2}}{v_{0}v_{1}}-\dfrac{1}{v_{0}}-\dfrac{1}{v_{1}}+1$$
As  $w_{0}=w_{0}'$ and $w_{1}=w_{1}'$, then
$$H_{3}(L_{g},\mathbb{Z})=H_{3}(L_{f},\mathbb{Z}) = (\mathbb{Z}_{m_{3}})^{\alpha(\textbf{w})+1}.$$
Finally, using (4.0.5) we obtain $\mu(L_{g})=\mu(L_{f}).$
\hfill$\square$
\medskip

Now, we have the following theorem:

\begin{theorem} Given the data $(d, \mathbf{w}=(w_{0},w_{1},w_{2},w_{3},w_{4}))$ such that   $d=w_{0}+w_{1}+w_{2}+w_{3}+w_{4}-1.$ Consider the polynomial of chain-cycle type  
$$f=z_{0}^{a_{0}}+z_{0}z_{1}^{a_{1}}+z_{4}z_{2}^{a_{2}}+z_{2}z_{3}^{a_{3}}+z_{3}z_{4}^{a_{4}}$$ of degree $d,$  
that cuts out a K\"ahler-Einstein orbifold hypersurface $X_f$ in weighted projective space $\mathbb{P}(\mathbf{w})$ such that the weight vector satisfies 
$\mathbf{w}=(w_{0},w_{1},w_{2},w_{3},w_{4})=(m_{3}v_{0},m_{3}v_{1},m_{2}v_{2},m_{2}v_{3},m_{2}v_{4})$ and $\gcd(m_{2},m_{3})=1$ and $m_2m_3=d.$ Suppose the corresponding Sasakian-Einstein link $L_f$ is a rational homology 7-sphere. We have 

\begin{enumerate}
\item If $L_f$ admits a twin which also admits a Sasaki-Einstein metric, then  the link $L_{f^{T}}$ associated to  the  BH-transpose rule of $f$  is a rational homology sphere admitting a Sasaki-Einstein metric.
\item If $L_f$ is such that  two elements of $\{ a_{2},a_{3},a_{4}\}$ are equal, then  the link $L_{f^{T}}$ associated to  the  BH-transpose rule of $f$  is a rational homology sphere admitting a Sasaki-Einstein metric.
\end{enumerate}

\end{theorem}

\noindent{\it Proof.} First, let us prove $(1).$ Notice that if the link associated to chain-cycle type singularity admits a twin, then the collection $a_{2}$,$a_{3},$ $a_{4}$ must be  different in pairs. Let  $\textbf{w}$ and $\textbf{w}'$ the weights associated to the  polynomials:
$$f=z_{0}^{a_{0}}+z_{0}z_{1}^{a_{1}}+z_{4}z_{2}^{a_{2}}+z_{2}z_{3}^{a_{3}}+z_{3}z_{4}^{a_{4}} \mbox{  and  } g=z_{0}^{a_{0}}+z_{0}z_{1}^{a_{1}}+z_{4}z_{2}^{a_{2}}+z_{2}z_{3}^{a_{4}}+z_{3}z_{4}^{a_{3}},$$
respectively such that the corresponding links $L_f$ and $L_g$ are twins. Then their associated matrices are
\begin{equation*}
    A_{f} =\begin{bmatrix}
      a_{0} & 0 & 0 & 0 & 0 \\
      1 & a_{1} & 0 & 0 & 0 \\
      0 & 0 & a_{2} & 0 & 1 \\
      0 & 0 & 1 & a_{3} & 0 \\
      0 &  0 & 0 & 1 & a_{4} 
    \end{bmatrix}, \ \ 
    A_{g} =\begin{bmatrix}
      a_{0} & 0 & 0 & 0 & 0 \\
      1 & a_{1} & 0 & 0 & 0 \\
      0 & 0 & a_{2} & 0 & 1 \\
      0 & 0 & 1 & a_{4} & 0 \\
      0 &  0 & 0 & 1 & a_{3} 
    \end{bmatrix}
\end{equation*}
Using the method of Berglund- Hübsch we obtain their corresponding transposes 
\begin{equation*}
    A_{f}^{T}=\begin{bmatrix}
      a_{0} & 1 & 0 & 0 & 0 \\
      0 & a_{1} & 0 & 0 & 0 \\
      0 & 0 & a_{2} & 1 & 0 \\
      0 & 0 & 0 & a_{3} & 1 \\
      0 &  0 & 1 & 0 & a_{4} 
    \end{bmatrix}, \ \
    A_{g}^{T}=\begin{bmatrix}
      a_{0} & 1 & 0 & 0 & 0 \\
      0 & a_{1} & 0 & 0 & 0 \\
      0 & 0 & a_{2} & 1 & 0 \\
      0 & 0 & 0 & a_{4} & 1 \\
      0 &  0 & 1 & 0 & a_{3} 
    \end{bmatrix}
\end{equation*}
and their associated polynomials $$\Tilde{f}^{T}=z_{0}^{a_{0}}z_{1}+z_{1}^{a_{1}}+z_{3}z_{2}^{a_{2}}+z_{4}z_{3}^{a_{3}}+z_{2}z_{4}^{a_{4}}$$ and $$\Tilde{g}^{T}=z_{0}^{a_{0}}z_{1}+z_{1}^{a_{1}}+z_{3}z_{2}^{a_{2}}+z_{4}z_{3}^{a_{4}}+z_{2}z_{4}^{a_{3}}.$$ We claim that if the link $L_{f}$ obtained by $f$ admits a Sasaki-Einstein, then the link $L_{g^{T}}$ admits a Sasaki-Einstein metric as well. So applying the argument once more, this time to $L_g$, we can conclude that $L_{f^{T}}$ admits a Sasaki-Einstein metric.

First, by the calculations performed in the previous section,  the degree $\Tilde{d}$ and the weight vector $\textbf{w}_{g^{T}}$ of $A_{g}^{T}$ are given by:
$$\textbf{w}_{g^{T}}=(\Tilde{w}_{0},\Tilde{w}_{1},\Tilde{w}_{2},\Tilde{w}_{3},\Tilde{w}_{4}) \ \ \ \mbox{and } \ \ \Tilde{d}=m_{3}m_{2}(m_{2}-1)$$
where
$$\Tilde{w}_{0}=m_{3}v_{1}(a_{1}-1), \Tilde{w}_{1}=m_{3}m_{2}v_{1}$$
and
\begin{align*}
    \Tilde{w}_{2} & = m_{2}(m_{2}-1)(1-a_{3}+a_{3}a_{4}) = (m_{2}-1)w_{2},\\
    \Tilde{w}_{3} & = m_{2}(m_{2}-1)(1-a_{4}+a_{2}a_{4}) = (m_{2}-1)w_{3},\\
    \Tilde{w}_{4} & = m_{2}(m_{2}-1)(1-a_{2}+a_{2}a_{3}) = (m_{2}-1)w_{4}.
\end{align*}
Notice that the index for $\textbf{w}_{g^{T}}$ is $\Tilde{I}=m_{2}-1$.  Let us  verify the inequality given in Theorem 2.1
\begin{equation}
\Tilde{I}\Tilde{d}< \frac{4}{3}\min\{ \Tilde{w}_{i}\Tilde{w}_{j} \}.
\end{equation}

As $m_{2}=a_{1}v_{1}+1$, then $\Tilde{w}_{0}<\Tilde{w}_{1}$. Moreover, we can assume, without loss of generality, that $\Tilde{w}_{2}\leq \Tilde{w}_{3}\leq \Tilde{w}_{4}$. Then the number $\min\{ \Tilde{w}_{i}\Tilde{w}_{j}\}$ will be equal to any of the following $\Tilde{w}_{0}\Tilde{w}_{1}$ or $\Tilde{w}_{0}\Tilde{w}_{2}$ or $\Tilde{w}_{2}\Tilde{w}_{3}$. Let us see that on each case (4.0.10) holds.

\begin{enumerate}
    \item[a)]  First, let us suppose $\min\{ \Tilde{w}_{i}\Tilde{w}_{j}\}=\Tilde{w}_{0}\Tilde{w}_{1}$. Then we have
    $\Tilde{I}\Tilde{d}=m_{3}m_{2}(m_{2}-1)^{2}$ and  
 $\Tilde{w}_{0}\Tilde{w}_{1}=m_{3}^{2}m_{2}v_{1}^{2}(a_{1}-1).$
Replacing in (4.0.10), we have 
$$m_{3}m_{2}(m_{2}-1)^{2} < \dfrac{4}{3}m_{3}^{2}m_{2}v_{1}^{2}(a_{1}-1).$$
Simplifying and using $v_{1}=\frac{m_{2}-1}{a_{1}}$, we obtain that the inequality in (4.0.10) is equivalent to the inequality 
\begin{equation}
1<\dfrac{4m_{3}(a_{1}-1)}{3a_{1}^{
2}}.
\end{equation}
The list of Johnson and Koll\'ar in \cite{JK}  produced  K\"ahler-Einstein 3-folds of index 1, such that all weight vectors, except the one  with weights $(13,143,775,620,465)$,  satisfy the inequality $d<w_{0}w_{1}.$  Since the weight vector $(13,143,775,620,465)$ produces the  chain-cycle polynomial $z_0^{155}+z_0z_1^{14}+z_4z_2^2+z_2z_3^2+z_3z_4^3$ with  $a_2=a_3.$ Thus  we can assume $d<w_{0}w_{1}$ for all the weight vectors under discussion.
  As $d=m_{3}m_{2}$, $w_{0}=m_{3}$ and $w_{1}=m_{3}v_{1}$,  we obtain 
$m_{3}m_{2}<(m_{3})^{2}v_{1}$ if and only if $m_{2}<m_{3}v_{1}.$
Putting $m_{2}=a_{1}v_{1}+1$ and replacing above, we obtain $1<(m_{3}-a_{1})v_{1}$, which implies $a_{1}<m_{3}$. Now, if we take $a_{1}\geq 4$, we obtain the inequalities 
$1\leq \dfrac{4(a_{1}-1)}{3a_{1}}$ and $1<\dfrac{m_{3}}{a_{1}}.$  Thus $1<\dfrac{4m_{3}(a_{1}-1)}{3a_{1}^{
2}}.$
On the other hand, if $a_{1}=2$ or $a_{1}=3$, from  (4.0.11) we have
 \begin{itemize}
 \item $a_{1}=2$ then  $1<\dfrac{m_{3}}{3}$ or equivalently $3<m_{3}.$
\item $a_{1}=3$ then  $1<\dfrac{8m_{3}}{27}$ or equivalently $\dfrac{27}{8}<m_{3}.$
\end{itemize}
Both inequalities hold for all weight vector.  Therefore, in this case, the link $L_{g^{T}}$ admits a Sasaki-Einstein metric.

\item[b)] Suppose $\min\{ \Tilde{w}_{i}\Tilde{w}_{j}\}=\Tilde{w}_{0}\Tilde{w}_{2}$.  Here we have the following equalities 
    $\Tilde{I}\Tilde{d}=m_{3}m_{2}(m_{2}-1)^{2}$ and 
 $\Tilde{w}_{0}\Tilde{w}_{2}=m_{3}v_{1}(a_{1}-1)(m_{2}-1)w_{2}.$
As $w_{2}=m_{2}v_{2}$, we obtain that 
    $\Tilde{I}\Tilde{d}<\dfrac{4}{3}\Tilde{w}_{0}\Tilde{w}_{2}$ which is equivalent to  $m_{3}m_{2}(m_{2}-1)^{2} <\dfrac{4}{3}m_{3}v_{1}(a_{1}-1)(m_{2}-1)m_{2}v_{2}.$
 Since $m_{2}-1=a_{1}v_{1}$, we simplify the above inequality:
$a_{1}<\dfrac{4}{3}(a_{1}-1)v_{2}.$
As $\dfrac{3a_{1}}{4(a_{1}-1)}<2$ and $ v_{2}\geq 2$, we can verify the above inequality.

\item[c)] Suppose  Suppose $\min\{ \Tilde{w}_{i}\Tilde{w}_{j}\}=\Tilde{w}_{2}\Tilde{w}_{3}$.  Here we have the equalities 
    $\Tilde{I}\Tilde{d}=m_{3}m_{2}(m_{2}-1)^{2} =d(m_{2}-1)^{2}$  and 
 $\Tilde{w}_{2}\Tilde{w}_{3}=(m_{2}-1)^{2}w_{2}w_{3}.$ Then
 $\Tilde{I}\Tilde{d}<\dfrac{4}{3}\Tilde{w}_{2}\Tilde{w}_{3}$ which is equivalent to $d<\dfrac{4}{3}w_{2}w_{3}.$
 This last inequality holds since the metric of link $L_{f}$ is Einstein.
\end{enumerate}

Now, let us prove $(2)$: if  two elements of $\{ a_{2},a_{3},a_{4}\}$ are equal,  $\textbf{w}$ does not admit a twin. We can suppose without loss of generality that $a_{2}=a_{3}$. Then the associated exponent matrix $A_{f}$ and its transpose are:

\begin{equation*}
    A_{f} =\begin{bmatrix}
      a_{0} & 0 & 0 & 0 & 0 \\
      1 & a_{1} & 0 & 0 & 0 \\
      0 & 0 & a_{2} & 0 & 1 \\
      0 & 0 & 1 & a_{2} & 0 \\
      0 &  0 & 0 & 1 & a_{4} 
    \end{bmatrix}, \ \ 
    A_{f}^{T} =\begin{bmatrix}
      a_{0} & 1 & 0 & 0 & 0 \\
      0 & a_{1} & 0 & 0 & 0 \\
      0 & 0 & a_{2} & 1 & 0 \\
      0 & 0 & 0 & a_{2} & 1 \\
      0 &  0 & 1 & 0 & a_{4} 
    \end{bmatrix}
\end{equation*}
We remember that the weights $w_{0},w_{1},\dots,w_{4}$ are expressed as
\begin{align*}
    w_{0} & = m_{3},\\
    w_{1} & = m_{3}v_{1},\\
    w_{2} & = m_{2}(1-a_{2}+a_{2}a_{4}),\\
    w_{3} & = m_{2}(1-a_{4}+a_{2}a_{4}),\\
    w_{4} & = m_{2}(1-a_{2}+a_{2}^{2}).
\end{align*}
and the weights $\Tilde{w}_{0},\Tilde{w}_{1},\dots\Tilde{w}_{4}$ are 
\begin{align*}
    \Tilde{w}_{0} & = m_{3}v_{1}(a_{1}-1),\\
    \Tilde{w}_{1} & = m_{3}m_{2}v_{1},\\
    \Tilde{w}_{2} & = m_{2}(1-a_{4}+a_{2}a_{4})=(m_{2}-1)w_{3},\\
    \Tilde{w}_{3} & = m_{2}(1-a_{2}+a_{2}a_{4})=(m_{2}-1)w_{2},\\
    \Tilde{w}_{4} & = m_{2}(1-a_{2}+a_{2}^{2})=(m_{2}-1)w_{4}.
\end{align*}

Since $\Tilde{w}_{0},\dots,\Tilde{w}_{4}$ are expressed in a similar way that in the previous case, we can apply the same process for almost all weight vectors $\textbf{w}$. The only exception is  the vector $\textbf{w}=(13,143,775,620,465).$ In this case, we calculate the weight vector $\textbf{w}_{f^{T}}=(169,2015,8680,10850,6510)$. Clearly, it verifies the formula (4.0.6).
\hfill$\square$

\begin{remark} From the two lists of 184 and 52 Sasaki-Einstein rational homology spheres found in \cite{BGN2} and in \cite{CL} respectively, we found  75 are given by polynomials of chain-cycle type: 29 different twins and 17 that do not posses twins, thus, according to Theorem 4.2, we obtain a total 58+17=75 new rational homology 7-spheres admitting Sasaki-Einstein metrics. All these new examples which are  links of not well-formed hypersurface singularites, are listed in a table given in the appendix.  
\end{remark}

In the next example we show how to produce these new Sasaki-Einstein rational homology 7-spheres in the context of the previous remark. 

\begin{example} Let us  consider the following data $(\mathbf{w}_f=(929,1858,2849,63,805), d_f=6503)$ found in Table 1 in \cite{CL}. These data  determine a  K\"ahler-Einstein 3-fold of index 1. Its associated polynomial is of  chain-cycle type: $f=z_0^7+z_0z_1^3+z_4z_2^2+z_2z_3^{58}+z_3z_4^8.$ One can find a twin for this element,  the data for it is given by  $(\mathbf{w}_g=(929,1858,3199,413,105), d_g=6503)$ which has an associated polynomial given by the following chain-cycle type singularity  $g=z_0^7+z_0z_1^3+z_4z_2^2+z_2z_3^8+z_3z_4^{58}.$  Using the BH-transpose rule on $f$ and $g$,  and after some calculations we obtain: 
\begin{itemize}
\item For $f$: the not well-formed data $(\mathbf{w}_{f^T}=(1858, 6503,9597, 315,1239), d_{f^T}=19509)$ with hypersurface singularity
$$f^T=z_{0}^{3}+z_{0}z_{1}^{7}+z_{4}z_{2}^{8}+z_{2}z_{3}^{58}+z_{3}z_{4}^{2}.$$ 

\item For $g$: the not well-formed data $(\mathbf{w}_{g^T}=(1858, 6503,8547, 2415,189), d_{g^T}=19509)$ with hypersurface singularity 
$$g^T=z_{0}^{3}+z_{0}z_{1}^{7}+z_{4}z_{2}^{58}+z_{2}z_{3}^{8}+z_{3}z_{4}^{2}.$$
\end{itemize}
It is not difficult to verify, thanks to  Orlik's algorithm, that the links associated to these invertible polynomials are rational homology 7-spheres that are twins: $H_3(L, \mathbb Z)=\mathbb Z_{929}$, and  Milnor number $\mu=17632$. Using the inequality given in Theorem 2.1 we conclude that both of them admit Sasaki-Einstein metrics. Thus we obtain these new Sasaki-Einstein rational homology 7-spheres.
\end{example}


The following two examples, point out that certain links that are not rational homology spheres can produce rational homology spheres admitting Sasaki-Einstein metrics via the BH-transpose rule. 
\begin{example} 
Consider the weight vector $\textbf{w}=(15,35,15,9,32)$ of the Johnson-Koll\'ar list. We can find a  BP-chain type polynomial of degree 105:
$$f=z_{0}^{7}+z_{1}^{3}+z_{2}^{7}+z_{2}z_{3}^{10}+z_{3}z_{4}^{3}.$$ It is not difficult to check that the link is of the form $24\#(S^3\times S^4)$ and it admits Sasaki-Einstein metric, see \cite{CL}.  
Its associated matrix of exponents is
\begin{equation*}
    A =\begin{bmatrix}
      7 & 0 & 0 & 0 & 0 \\
      0 & 3 & 0 & 0 & 0 \\
      0 & 0 & 7 & 0 & 0 \\
      0 & 0 & 1 & 10 & 0 \\
      0 &  0 & 0 & 1 & 3 
    \end{bmatrix}
\end{equation*}
Through the BH transpose rule we find dual data $\tilde{\textbf{w}} = (15,35,14,7,35)$ with $\tilde{d}=105$. Using Orlik's formula, we found that its link $L_{f^T}$ is a rational homology 7-sphere over a not well-formed orbifold: its homology and Milnor number are $H_{3}(L, \mathbb Z)=(\mathbb{Z}_{7})^{26}$ and $\mu=2184$ respectively. Furthermore, from the inequality given in Theorem 2.1, we conclude that this link admits Sasaki-Einstein metric.

\end{example}

\begin{example}
    Consider the vector weight $\textbf{w}=(5,35,57,64,160)$ and the degree $d=320$ of the list of Johnson and Koll\'ar. It    
    is not known if this element is tiger or  K\"ahler-Einstein, hence it is unknown whether 
    the corresponding link, which is of  the form $36\#(S^3\times S^4)$, admits or not Sasaki-Einstein metric.  Rearranging the  weight vector $\textbf{w}=(64,160,5,35,57)$ we have the BP-chain polynomial $f$ of degree $d=320$:
$$f=z_{0}^{5}+z_{1}^{2}+z_{1}z_{2}^{32}+z_{2}z_{3}^{9}+z_{3}z_{4}^{5}$$
Its associated matrix of exponent is
\begin{equation*}
    A =\begin{bmatrix}
      5 & 0 & 0 & 0 & 0 \\
      0 & 2 & 0 & 0 & 0 \\
      0 & 1 & 32 & 0 & 0 \\
      0 & 0 & 1 & 9 & 0 \\
      0 &  0 & 0 & 1 & 5 
    \end{bmatrix}
\end{equation*}

If we solve the matrix equation $A^{t}\tilde{W} = \tilde{D}$, we obtain the  weight vector $\tilde{\textbf{w}}=(576,1399,82,256,576)$ with $\tilde{d}=2880$. From  Orlik's formula we obtain  $H_{3}=\mathbb{Z}_{90}\oplus(\mathbb{Z}_{18})^{3}$ 
 and $\mu=5924$. In other words, we obtain a rational homology sphere. Furthermore, by the inequality of Theorem 2.1, it follows  that the link $L$ associated to $\tilde{w}$ is Sasaki-Einstein.

\end{example}

The next example shows that the Sasaki-Einstein inequality given in Theorem 2.1 is not preserved for general 7-manifolds unless these are rational homology 7-spheres (compare with Example 2 in \cite{Go}).  

\begin{example} Consider the vector weight $\textbf{w}=(9,15,30,47,50)$ of the list of Johnson and Koll\'ar. It is easy to establish the existence of a Sasaki-Einstein metric on the link $L$ associated to it. Also, one obtains 
$H_3(L, \mathbb Z)=\mathbb{Z}^{24}$, thus $L$ is of the form $24\# (S^3\times S^4).$ Rearranging the vector weight $\textbf{w}=(30,50,15,9,47)$ we obtain a BP-chain polynomial $f=z_{0}^{5}+z_{1}^{3}+z_{2}^{10}+z_{2}z_{3}^{15}+z_{3}z_{4}^{3}$ of degree $d=150$ for these weights. Its associated matrix of exponent is given by 
\begin{equation*}
    A =\begin{bmatrix}
      5 & 0 & 0 & 0 & 0 \\
      0 & 3 & 0 & 0 & 0 \\
      0 & 0 & 10 & 0 & 0 \\
      0 & 0 & 1 & 15 & 0 \\
      0 & 0 & 0 & 1 & 3
    \end{bmatrix}
\end{equation*}

From the matrix equation $A^{t}\tilde{W} = \tilde{D}$, we obtain the vector weight $\tilde{\textbf{w}} = (90,150,43,20,150)$ with $\tilde{d}=450$. Using Orlik's formula we obtain  
$H_{3}(L_{f^T}, \mathbb Z)=(\mathbb{Z}_{50})^{2}\oplus(\mathbb{Z}_{10})^{6}.$ Thus, the  link $L_{f^T}$ is a rational homology sphere. In addition, we can verify that the link $L$ associated to $\tilde{\textbf{w}}$ does not satisfy the inequality of Theorem 2.1. So it is uncertain whether this link admits a Sasaki-Einstein metric or not.
\end{example}

The next example points out  that the BH-transpose rule might not produce data that satisfy the Sasaki-Einstein inequality given in Theorem 2.1.

\begin{example}
From the list of Johnson and Koll\'ar, consider the data  $(\mathbf{w}=(13,13,125,100,75), d=325)$ which produces  an index one 3-fold. It is not known whether this orbifold is tiger or  K\"ahler-Einstein, so the corresponding link $L$ does not necessarily admit an Einstein metric. We notice that there exist polynomials of type: BP-cycle, chain-cycle and cycle-cycle for these weights:
$$\mbox{BP-cycle}: \ \ f=z_{0}^{25}+z_{1}^{25}+z_{4}z_{2}^{2}+z_{2}z_{3}^{2}+z_{3}z_{4}^{3}$$
$$\mbox{chain-cycle}: \ \ f=z_{0}^{25}+z_{0}z_{1}^{24}+z_{4}z_{2}^{2}+z_{2}z_{3}^{2}+z_{3}z_{4}^{3}$$
$$\mbox{cycle-cycle}: \ \ f=z_{1}z_{0}^{24}+z_{0}z_{1}^{24}+z_{4}z_{2}^{2}+z_{2}z_{3}^{2}+z_{3}z_{4}^{3}$$
From  Orlik's formula we obtain $H_{3}(L,\mathbb{Z})=(\mathbb{Z}_{13})^{24}$. Thus, the link is a rational homology 7-sphere.
Next, we will use the BH-transpose rule  for the chain-cycle singularity. We calculate its matrix of exponents:
\begin{equation*}
    A=\begin{bmatrix}
        25 & 0 & 0 & 0 & 0 \\
        1 & 24 & 0 & 0 & 0 \\
        0 & 0 & 2 & 0 & 1 \\
        0 & 0 & 1 & 2 & 0 \\
        0 & 0 & 0 & 1 & 3
    \end{bmatrix}
\end{equation*}
Via the BH-transpose rule we obtain the dual data $\tilde{\textbf{w}}=(299,325,2400,3000,1800)$ and $\tilde{d}=7800.$ From  Orlik's formula for the corresponding link $L_{f^T}$  we obtain  $H_{3}(L_{f^T},\mathbb Z)=\mathbb{Z}_{13}.$ In this case, we can not determine if this link admits a Sasaki-Einstein metric since  these weights do not verify the inequality
$\Tilde{I}\Tilde{d}< \frac{4}{3}\min\{ \Tilde{w}_{i}\Tilde{w}_{j} \}.$
\end{example}

\begin{remark} In the context of   Corollary 3.4, for $f$  of cycle,  cycle-cycle or BP-cycle type such that the hypersurface $X_f$ cut out by it is a rational homology of a complex projective space, as argued in the introduction, one must obtain two rational homology 7-spheres $L_f$ and $L_{f^T}$ that, moreover,  are twins. Since the transverse spaces of these links are related by the HB-transpose rule, it is natural to look for more relations between them. For instance,  it is not difficult to see that, both orbifolds are  related to a Fermat hypersurface via Shioda maps \cite{Ke}. Indeed, consider 
  $f=\sum_{j=0}^4 \prod_{i=0}^4 x_i^{a_{i j}}$ be an invertible polynomial in $\mathbb C^5$
with defining  matrix $A=\left(a_{i j}\right)_{i, j}$ and denote by $X_f\subset \mathbb{P}\left(w_0, \ldots, w_4\right)$  the weighted hypersurface  of degree $d$ cut out by $f$ and by $X_{f^T}\subset \mathbb{P}\left(\tilde{w}_0, \ldots, \tilde{w}_4\right)$ the weighted hypersurface  of degree $\tilde{d}$ cut out by 
$f^T.$  Now consider the matrix $B=dA^{-1}=\left(b_{j k}\right)_{j, k}$ which has only integer entries. The Shioda maps are the rational maps 
$$
\begin{aligned}
\phi_A & : \mathbb{P}^4 \dashrightarrow  \mathbb{P}\left(w_0, \ldots, w_4\right) \\
\phi_{A^T} & : \mathbb{P}^4 \dashrightarrow \mathbb{P}\left(\tilde{w}_0, \ldots, \tilde{w}_4\right),
\end{aligned}
$$
where
$$
\begin{aligned}
& \left(y_0: \ldots: y_4\right) \stackrel{\phi_A}{\longmapsto}\left(x_0: \ldots: x_4\right), \quad x_j=\prod_{k=0}^4 y_k^{b_{j k}}, \\
& \left(y_0: \ldots: y_4\right) \stackrel{\phi_{A^T}}{\longmapsto}\left(z_0: \ldots: z_4\right), \quad z_j=\prod_{k=0}^4 y_k^{b_{k j}} .
\end{aligned}
$$
These maps restricted to the Fermat hypersurface $X_{d}\subset \mathbb{P}^4$ cut out by  the polynomial
$$
F_{d}:=\sum_{i=0}^4 y_i^d,
$$
give the rational maps $\phi_{A}: X_{d}\dashrightarrow X_f$ and $\phi_{A^T}: X_{d}\dashrightarrow X_{f^T}.$  It is well-known that  link $L_d$ associated to the Fermat hypersurface $X_d$,  admits a regular Sasakian geometry which is not Einstein in general: in order for that to happen, first of all,  the Fermat  hypersurface $F_{d}:=\sum_{i=0}^4 y_i^d$ must be Fano and this only occurs for $d<5$, condition that is not satisfied for the examples exhibited here. On the other hand if $d\geq 5,$ $F_{d}$ admits either null or negative K\"ahler-Einstein metric by Aubin-Yau theorem. It follows that the link admits regular Sasaki $\eta$-Einstein metrics which are negative or null \cite{BGM}. Thus, by  lifting the  Shioda maps to the associated links, one obtains   
 a regular negative Sasaki $\eta$-Einstein model related  to the quasi-regular Sasaki-Einstein structures on the twins.

\end{remark}

\section{Appendix}
In this appendix  we include a table that presents the  75 new examples of Sasaki-Einstein rational homology 7-spheres referred in Remark 4.1. This table lists the weights of the original data, the resulting weights produced by the BH-transpose rule, the degree, the Milnor number and finally the third homology groups of the corresponding links. Unlike the tables given in \cite{BGN2} and \cite{CL} most of the elements (58 to be precise) of the list have even degree and all of them are $S^1$-orbibundles over not well-formed orbifolds. We also give the following links to codes implemented in Matlab in order to avoid tedious computations: 
\begin{itemize}
\item Code4alt.m \href{https://github.com/Jcuadrosvalle/Codes-in-Matlab/blob/main/Code4alt.m}{https://github.com/Jcuadrosvalle/Codes-in-Matlab/blob/main/Code4alt.m} computes the degree, Milnor number and third homology group of the link of a hypersurface singularity of index greater or equal than 1.
\item Grado.m  \href{https://github.com/Jcuadrosvalle/Codes-in-Matlab/blob/main/Degree.m}{https://github.com/Jcuadrosvalle/Codes-in-Matlab/blob/main/Degree.m} computes the weights and degree of the polynomial from the entries of the matrix of exponents.
\end{itemize}
\medskip

\centering

\begin{center}
\noindent{\bf Table: Rational homology 7-spheres admitting S-E structures referred in Remark 4.1} 
\end{center}

 \begin{longtable}{| c || c | c | c | c | } 
\hline 
$\textbf{w}=(w_{0},w_{1},w_{2},w_{3},w_{4})$ & $\tilde{\textbf{w}}_{T}=(\tilde{w}_{0},\tilde{w}_{1},\tilde{w}_{2},\tilde{w}_{3},\tilde{w}_{4})$ &   $\tilde{d}$   &  $\tilde{\mu}$  & $H_{3}(L_{T},\mathbb{Z})$   \\  \hline 
\hline 
\endhead

 
$(73,73,95,45,80)$ & $(219,365,420,200,260)$ & $1460$ & $1224$ & $\mathbb{Z}_{73}$ \\ \hline

$(73,73,105,65,50)$ & $(219,365,380,320,180)$ & $1460$ & 1224 & $\mathbb{Z}_{73}$ \\ \hline

$(17,34,175,125,75)$ & $(187,425,1500,2100,900)$ & $5100$ & $4624$ & $\mathbb{Z}_{17}$ \\ \hline

$(67,67,161,28,147)$ & $(335,469,1302,210,504)$ & $2814$ & $2442$ & $\mathbb{Z}_{67}$ \\ \hline

$(67,67,217,84,35)$ & $(335,469,966,882,168)$ & $2814$ & $2442$ & $\mathbb{Z}_{67}$ \\ \hline

$(19,57,175,100,125)$ & $(133,475,1400,1000,800)$ & $3800$ & $3474$ & $\mathbb{Z}_{19}$ \\ \hline

$(118,118,185,135,35)$ & $(177,295,270,370,70)$ & $1180$ & $1989$ & $\mathbb{Z}_{590}$ \\ \hline

$(77,77,333,180,27)$ & $(539,693,1440,2664,216)$ & $5544$ & $4940$ & $\mathbb{Z}_{77}$ \\ \hline

$(253,253,545,40,175)$ & $(759,1265,2400,260,380)$  & $5060$ & $4284$ & $\mathbb{Z}_{253}$ \\ \hline

$(253,253,600,95,65)$ & $(759,1265,2180,700,160)$ & $5060$ & $4284$ & $\mathbb{Z}_{253}$ \\ \hline

$(113,226,715,377,39)$ & $(565,1469,2262,4290,234)$ & $8814$ & $8176$ & $\mathbb{Z}_{113}$ \\ \hline

$(64,512,475,375,175)$ & $(128,1600,1125,1425,525)$ & $4800$ & $4599$ & $\mathbb{Z}_{64}$ \\ \hline

$(127,381,559,52,533)$ & $(381,1651,3172,260,1144)$ & $6604$ & $6174$ & $\mathbb{Z}_{127}$ \\ \hline

$(127,381,793,286,65)$ & $(381,1651,2236,2132,208)$ & $6604$ & $6174$ & $\mathbb{Z}_{127}$ \\ \hline

$(373,373,780,35,305)$ & $(1119,1865,3620,220,640)$ & $7460$ & $6324$ & $\mathbb{Z}_{373}$ \\ \hline

$(373,373,905,160,55)$ & $(1119,1865,3120,1220,140)$ & $7460$ & $6324$ & $\mathbb{Z}_{373}$ \\ \hline

$(13,143,775,620,465)$ & $(169,2015,8680,10850,6510)$ & $28210$ & $25884$ & $\mathbb{Z}_{13}$ \\ \hline

$(701,701,381,123,198)$ & $(701,2103,786,276,342)$ & $4206$ & $3500$ & $\mathbb{Z}_{701}$ \\ \hline

$(701,701,393,171,138)$ & $(701,2103,762,396,246)$ & $4206$ & $3500$ & $\mathbb{Z}_{701}$ \\ \hline

$(17,238,889,635,381)$ & $(136,2159,5715,8001,3429)$ & $19431$ & $18160$ & $\mathbb{Z}_{17}$ \\ \hline

$(334,668,763,525,49)$ & $(668,2338,1575,2289,147)$ & $7014$ & $6327$ & $\mathbb{Z}_{334}$ \\ \hline

$(881,881,465,99,318)$ & $(881,2643,1014,216,534)$ & $5286$ & $4400$ & $\mathbb{Z}_{881}$ \\ \hline

$(881,881,507,267,108)$ & $(881,2643,930,636,198)$ & $5286$ & $4400$ & $\mathbb{Z}_{881}$ \\ \hline

$(73,584,1435,779,123)$ & $(292,2993,3895,7175,615)$ & $14965$ & $14472$ & $\mathbb{Z}_{73}$ \\ \hline

$(65,650,1581,867,153)$ & $(52,663,867,1581,153)$ & $3315$ & $16064$ & $\mathbb{Z}_{3315}\oplus(\mathbb{Z}_{51})^{3}$ \\ \hline

$(481,962,1519,77,329)$ & $(962,3367,4851,399,525)$ & $10101$ & $9120$ & $\mathbb{Z}_{481}$ \\ \hline

$(481,962,1617,175,133)$ & $(962,3367,4557,987,231)$ & $10101$ & $9120$ & $\mathbb{Z}_{481}$ \\ \hline

$(185,740,1911,987,63)$ & $(148,777,987,1911,63)$ & $3885$ & $18584$ & $\mathbb{Z}_{3885}\oplus(\mathbb{Z}_{21})^{3}$ \\ \hline

$(1331,1331,687,87,558)$ & $(1331,3993,1560,186,918)$ & $7986$ & $6650$ & $\mathbb{Z}_{1331}$ \\ \hline

$(1331,1331,780,459,93)$ & $(1331,3993,1374,1116,174)$ & $7986$ & $6650$ & $\mathbb{Z}_{1331}$ \\ \hline

$(29,667,1807,1112,417)$ & $(145,4031,6672,10842,2502)$ & $24186$ & $23212$ & $\mathbb{Z}_{29}$ \\ \hline

$(1457,1457,1011,120,327)$ & $(1457,4371,2112,294,510)$ & $8742$ & $7280$ & $\mathbb{Z}_{1457}$ \\ \hline

$(1457,1457,1056,255,147)$ & $(1457,4371,2022,654,240)$ & $8742$ & $7280$ & $\mathbb{Z}_{1457}$ \\ \hline

$(409,2045,1320,187,539)$ & $(409,4499,2838,484,770)$ & $8998$ & $8568$ & $\mathbb{Z}_{409}$ \\ \hline

$(409,2045,1419,385,242)$ & $(409,4499,2640,1078,374)$ & $8998$ & $8568$ & $\mathbb{Z}_{409}$ \\ \hline

$(43,1333,1875,500,1625)$ & $(129,5375,9500,2500,4000)$ & $21500$ & $20874$ & $\mathbb{Z}_{43}$ \\ \hline

$(43,1333,2375,1000,625)$ & $(129,5375,7500,6500,2000)$ & $21500$ & $20874$ & $\mathbb{Z}_{43}$ \\ \hline

$(2069,2069,1413,102,555)$ & $(2069,6207,3042,246,852)$ & $12414$ & $10340$ & $\mathbb{Z}_{2069}$ \\ \hline

$(2069,2069,1521,426,123)$ & $(2069,6207,2826,1110,204)$ & $12414$ & $10340$ & $\mathbb{Z}_{2069}$ \\ \hline

$(929,1858,2849,63,805)$ & $(1858, 6503, 1239, 9597, 315)$ & $19509$ & $17632$ & $\mathbb{Z}_{929}$ \\ \hline

$(929,1858,3199,413,105)$ & $(1858, 6503, 189, 8547, 2415)$ & $19509$ & $17632$ & $\mathbb{Z}_{929}$ \\ \hline

$(289,2312,2725,125,1775)$ & $(578, 7225, 2775, 10575, 525)$ & $21675$ & $21024$ & $\mathbb{Z}_{289}$ \\ \hline

$(289,2312,3525,925,175)$ & $(578, 7225, 375, 8175, 5325)$ & $21675$ & $21024$ & $\mathbb{Z}_{289}$ \\ \hline

$(1297,3891,2653,119,1120)$ & $(1297, 9079, 5950, 308, 1526)$ & $18158$ & $16848$ & $\mathbb{Z}_{1297}$ \\ \hline

$(1297,3891,2975,763,154)$ & $(1297, 9079, 5306, 2240, 238)$ & $18158$ & $16848$ & $\mathbb{Z}_{1297}$ \\ \hline

$(217,4557,2752,731,1075)$ & $(217, 9331, 5590, 1892, 1634)$ & $18662$ & $18360$ & $\mathbb{Z}_{217}$ \\ \hline

$(217,4557,2795,817,946)$ & $(217, 9331, 5504, 2150, 1462)$ & $18662$ & $18360$ & $\mathbb{Z}_{217}$ \\ \hline

$(49,1862,4393,2483,573)$ & $(196, 9359, 12415, 21965, 2865)$ & $46795$ & $45648$ & $\mathbb{Z}_{49}$ \\ \hline

$(3401,3401,2298,93,1011)$ & $(3401, 10203, 5046, 222, 1536)$ & $20406$ & $17000$ & $\mathbb{Z}_{3401}$ \\ \hline

$(3401,3401,2523,768,111)$ & $(3401, 10203, 4596, 2022, 186)$ & $20406$ & $17000$ & $\mathbb{Z}_{3401}$ \\ \hline

$(129,3612,4165,425,2635)$ & $(86, 3655, 5185, 595, 1445)$ & $10965$ & $32384$ & $\mathbb{Z}_{10965}\oplus\mathbb{Z}_{85}$ \\ \hline

$(129,3612,5185,1445,595)$ & $(86, 3655, 4165, 2635, 425)$ & $10965$ & $32384$ & $\mathbb{Z}_{10965}\oplus\mathbb{Z}_{85}$ \\ \hline

$(657,3942,4693,95,3097)$ & $(438, 4161, 6175, 133, 1577)$ & $12483$ & $36080$ & $\mathbb{Z}_{12483}\oplus\mathbb{Z}_{19}$ \\ \hline

$(657,3942,6175,1577,133)$ & $(438, 4161, 4693, 3097, 95)$ & $12483$ & $36080$ & $\mathbb{Z}_{12483}\oplus\mathbb{Z}_{19}$ \\ \hline

$(1135,5675,3476,143,2057)$ & $(1135, 12485, 8206, 352, 2794)$ & $24970$ & $23814$ & $\mathbb{Z}_{1135}$ \\ \hline

$(1135,5675,4103,1397,176)$ & $(1135, 12485, 8206, 2794, 352)$ & $24970$ & $23814$ & $\mathbb{Z}_{1135}$ \\ \hline

$(1505,6020,3357,2547,117)$ & $(1505, 13545, 234, 5094, 6714)$ & $27090$ & $25568$ & $\mathbb{Z}_{1505}$ \\ \hline

$(3532,7064,5355,115,1595)$ & $(1766, 8830, 1075, 5835, 155)$ & $17660$ & $31779$ & $\mathbb{Z}_{17660}$ \\ \hline

$(3532,7064,5835,1075,155)$ & $(1766, 8830, 1595, 115, 5355)$ & $17660$ & $31779$ & $\mathbb{Z}_{17660}$ \\ \hline

$(141,9729,4031,2224,3475)$ & $(141, 19599, 8618, 4726, 6116)$ & $39198$ & $38780$ & $\mathbb{Z}_{141}$ \\ \hline

$(141,9729,4309,3058,2363)$ & $(141, 19599, 8062, 6950, 4448)$ & $39198$ & $38780$ & $\mathbb{Z}_{141}$ \\ \hline

$(113,8362,9589,1115,6021)$ & $(226, 25199, 35457, 4683, 10035)$ & $75597$ & $74704$ & $\mathbb{Z}_{113}$ \\ \hline

$(113,8362,11819,3345,1561)$ & $(226, 25199, 28767, 18063, 3345)$ & $75597$ & $74704$ & $\mathbb{Z}_{113}$ \\ \hline

$(1351,12159,6859,209,5092)$ & $(1351, 25669, 16948, 494, 6878)$ & $51338$ & $49950$ & $\mathbb{Z}_{1351}$ \\ \hline

$(1351,12159,8474,3439,247)$ & $(1351, 25669, 13718, 10184, 418)$ & $51338$ & $49950$ & $\mathbb{Z}_{1351}$ \\ \hline

$(177,15399,7175,5950,2275)$ & $(177, 30975, 11900, 14350, 4550)$ & $61950$ & $61424$ & $\mathbb{Z}_{177}$ \\ \hline

$(193,18335,10887,3247,4202)$ & $(193, 36863, 21774, 8404, 6494)$ & $73726$ & $73152$ & $\mathbb{Z}_{193}$ \\ \hline

$(2416,19328,10965,187,8177)$ & $(1208, 20536, 13617, 221, 5491)$ & $41072$ & $79695$ & $\mathbb{Z}_{41072}$ \\ \hline

$(2416,19328,13617,5491,221)$ & $(1208, 20536, 10965, 8177, 187)$ & $41072$ & $79695$ & $\mathbb{Z}_{41072}$ \\ \hline

$(217,23219,13115,2795,7310)$ & $(217, 46655, 28810, 6880, 10750)$ & $93310$ & $92664$ & $\mathbb{Z}_{217}$ \\ \hline

$(217,23219,14405,5375,3440)$ & $(217, 46655, 26230, 14620, 5590)$ & $93310$ & $92664$ & $\mathbb{Z}_{217}$ \\ \hline

$(316,24648,13345,1727,9577)$ & $(158, 24806, 15857, 2041, 6751)$ & $49612$ & $98595$ & $\mathbb{Z}_{49612}$ \\ \hline

$(316,24648,15857,6751,2041)$ & $(158, 24806, 13345, 9577, 1727)$ & $49612$ & $98595$ & $\mathbb{Z}_{49612}$ \\ \hline

$(301,44849,24219,3289,17342)$ & $(301, 89999, 57408, 7774, 24518)$ & $179998$ & $179100$ & $\mathbb{Z}_{301}$ \\ \hline

$(301,44849,28704,12259,3887)$ &  $(301, 89999, 48438, 34684, 6578)$ & $179998$ & $179100$ & $\mathbb{Z}_{301}$ \\ \hline
\end{longtable}

\begin{remark} As in  \cite{BGN3}, one can produce homotopy $9$-spheres bounding parallelizable manifolds that admit positive Ricci curvature from each of the rational homology 7-spheres in the table presented above: 
one  considers links  $L_g$ as  $p$-fold branched covers of $S^{2m+1}$ ramified over rational homology spheres, where $g$ is a polynomial of the form 
$$
g=z_0^p+f\left(z_1, \ldots, z_m\right)
$$
for $p>1$ and  where $f$ determines one of the links in the table given above. 
It is well-known that in dimension 9, the only diffeomorphism type that can occur in this scenario is either the standard 9-sphere or the Kervaire 9-sphere. Moreover, due to a  result of  Levine's (see [21], page 69) one can  determine the diffeomorphism type of these links via the Alexander polynomial: If   $\Delta_g(-1) \equiv \pm 1(\bmod 8)$ the link is diffeomorphic to the standard 9-sphere, and  if $\Delta_g(-1) \equiv \pm 3(\bmod 8)$ then $L_g$ is  diffeomorphic to the Kervaire $9$-spheres. Thus one can determine the diffeomorphism type in terms of the weight $p$. Furthermore, applying similar ideas as the ones given in Remark 1 in \cite{CL}, one can find bounds for the weight $p$ such that $L_g$ does not admit extremal Sasaki metrics. It is worth mentioning that a recent result in \cite{LST} establishes the existence of Sasaki-Einstein metrics on all odd-dimensional homotopy spheres which bound parallelizable manifolds and  the metrics found in that article are constructed on links over Brieskorn-Pham polynomials and not over the more general setting of invertible polynomials, which are the ones exhibited here. 
\end{remark}

\section*{Declarations}





\subsection*{Funding} The first author received  financial support from Pontificia Universidad Católica del Perú through project VRI-DFI 2016-1-0060.



\end{document}